\documentclass{lms}
\usepackage{amssymb}
\usepackage{amsmath}

\newcommand{\an}{\mbox{\scriptsize an}}
\newcommand{\Sing}{\mbox{\scriptsize Sing}}
\newcommand{\supp}{\mbox{\scriptsize supp}}
\newcommand{\fin}{\mbox{\scriptsize fin}}
\newcommand{\Arg}{\mbox{\scriptsize Arg}}
\newcommand{\dist}{\mbox{\scriptsize dist}}

\newtheorem{theorem}{Theorem}[section] 

\newnumbered{assertion}{Assertion}    
\newnumbered{conjecture}{Conjecture}  
\newnumbered{definition}{Definition}
\newnumbered{hypothesis}{Hypothesis}
\newnumbered{remark}{Remark}
\newnumbered{note}{Note}
\newnumbered{observation}{Observation}
\newnumbered{problem}{Problem}
\newnumbered{question}{Question}
\newnumbered{algorithm}{Algorithm}
\newnumbered{example}{Example}
\newunnumbered{notation}{Notation} 

\title[The Riemann Mapping Theorem for semianalytic domains and o-minimality
]
 {The Riemann Mapping Theorem for semianalytic domains and o-minimality} 

\author{T. Kaiser}

\classno{03C64, 32B20, 30C20, 30E15, 30D60, 30D05, 37E35}

\extraline{{\it Keywords and phrases}: o-minimal structures, semianalytic sets,
Riemann Mapping Theorem, quasianalytic classes, dynamical systems\\
Supported by DFG-project KN202/5-2}

\begin{document}
\maketitle

\begin{abstract}
We consider the Riemann Mapping Theorem in the case of a bounded
simply connected and semianalytic domain. We show that the germ at 0
of the Riemann map (i.e. biholomorphic map) from the upper half
plane to such a domain can be realized in a certain quasianalytic
class if the angle of the boundary at the point to which 0 is
mapped, is greater than 0. This quasianalytic class was introduced
and used by Ilyashenko in his work on Hilbert's 16$^{\mbox{\tiny
th}}$ problem. With this result we can prove that the Riemann map
from a bounded simply connected semianalytic domain onto the unit
ball is definable in an o-minimal structure, provided that at
singular boundary points the angles of the boundary are irrational
multiples of $\pi$.
\end{abstract}

\vspace{1cm}\noindent {\large \bf Introduction}

\medskip\noindent
One of the central theorems in complex analysis is the

\vspace{0.5cm}\noindent {\bf Riemann Mapping Theorem}

\medskip\noindent
{\it Let $\Omega \subsetneqq {\mathbb C}$ be a simply connected
domain in the plane. Then $\Omega$ can be mapped biholomorphically
onto the unit ball} $B(0, 1)$.

\medskip\noindent
A nice overview of its proofs and their history can be found in
Remmert [29].

\noindent `Riemann maps' are in general transcendental and not
algebraic functions. The goal of this paper is to show a `tame'
content of the Riemann Mapping Theorem. An important framework of
`tame' geometry is given by the category of semialgebraic sets and
functions. Semialgebraic sets and functions have very nice
finiteness properties (see Bochnak et al. [3]). A more general
framework to study sets with singularities which still exhibit a
nice behaviour is given by the category of subanalytic sets and
functions (see Bierstone-Milman [2], Denef-Van den Dries [6], \L
ojasiewicz [24] and Shiota [31]). But we lose some of the finiteness
properties we have in semialgebraic geometry. These are still valid,
however, if we restrict ourselves to globally subanalytic sets and
functions, i.e. to sets and functions which are subanalytic in the
ambient projective space (compare with Van den Dries [7] and Van den
Dries-Miller [10]). We lose little here since all bounded
subanalytic sets are globally subanalytic.

\medskip\noindent
The question is now the following:\\
Given a bounded simply connected domain in the plane which is
semianalytic (in ${\mathbb R}^2$ subanalytic is the same as
semianalytic; see the preliminary section), what can be said about
the Riemann map from the domain to the unit ball with respect to
'tame' behaviour. Efroymson showed in [15] that any semialgebraic
and simply connected domain in the plane is Nash isomorphic to
${\mathbb R}^2$. The Riemann map however, is in general not
subanalytic as we will see below. It often occurs that the category
of (globally) subanalytic sets and maps is too small for concepts
from analysis. For example, solutions to ordinary differential
equations with subanalytic row data and volume of subanalytic
families are in general not subanalytic. But in these two examples
the resulting functions can be realized in so-called o-minimal
structures (see Comte et al. [4] and Wilkie [33]). O-minimal
structures represent an excellent 'tame' generalization of the
category of semialgebraic or globally subanalytic sets and functions
and are defined by finiteness properties (see Van den Dries [8] and
the preliminaries for more details). Under certain conditions on the
singular boundary points of the given bounded, semianalytic domain
we obtain an o-minimality result:

\bigskip\noindent
{\bf Theorem A}

\noindent
{\it There is an o-minimal structure with the following property:\\
Let $\Omega \subset {\mathbb C}$ be a bounded simply connected
domain which is semianalytic. Suppose that the following condition
holds: if $x$ is a singular boundary point of $\Omega$ then the
angle of the boundary at $x$ is an irrational multiple of $\pi$.
Then the Riemann map from $\Omega$ onto the unit ball (i.e. its
graph considered as a subset of ${\mathbb R}^4$) is definable in
this o-minimal structure.}

\medskip\noindent
We will make the notion of angle at a boundary point more precise in
the text. As applications we obtain from Theorem A the following:
working with polygons we get the definability of the
Schwarz-Christoffel maps, and working with circular polygons we get
the definability of certain ratios of hypergeometric functions in
this o-minimal structure. Using general model theory the Riemann
Mapping Theorem can be transferred by Theorem A to real closed
fields of arbitrary large cardinality in a 'tame' way (compare with
Peterzil-Starchenko [26, 27] and also with Huber-Knebusch [17] and
Knebusch [22] for the development of complex analysis in
o-minimal structures on arbitrary real closed fields).\\
The ideas of the proof are the following: let $\Omega$ be a bounded
semianalytic domain in the plane and let $\varphi \colon \Omega \to
B(0, 1)$ be a Riemann map (i.e. a biholomorphic map). Let $x \in
\overline{\Omega}$. If $x \in \Omega$ then $\varphi$ is analytic in
a neighbourhood of $x$ and if $x \in
\partial \Omega$ is a nonsingular boundary point of $\Omega$
(i.e. the boundary is an analytic manifold at $x$) then $\varphi$
has an analytic extension to a neighbourhood of $x$ by the Schwarz
reflection principle. So the interesting (and hard) case is when $x
\in
\partial \Omega$ is singular. Taking the inverse of $\varphi$ and composing it
with a M\"obius transformation (a linear fractional map) it is an
equivalent problem to consider a Riemann map $\Phi \colon {\mathbb
H} \to \Omega$ where ${\mathbb H}$ denotes the upper half plane.

\noindent Given a simply connected domain $D$ which has an analytic
corner at $0 \in
\partial D$ (i.e. the boundary at 0 is given by two regular analytic arcs
which intersect in an angle $\sphericalangle D$ greater than 0),
Lehman showed in [23] (see also Pommerenke [28, p.58]) that a
Riemann map $\Phi \colon {\mathbb H} \to D$ with $\Phi (0) = 0$ has
an asymptotic development at 0 of the following kind:
$$\Phi (z) \sim \sum\limits_{n = 0}^{\infty} a_n P_n (\log z) z^{\alpha_n}
\quad \mbox{as} \quad z \longrightarrow 0 \quad \mbox{on} \quad
{\mathbb H}, \leqno(\dagger)$$

\noindent i.e. for each $N \in {\mathbb N}_0$ we have
$$\Phi (z) - \sum\limits_{n = 0}^N a_n P_n (\log z) z^{\alpha_n} = o
(z^{\alpha_N}) \quad \mbox{as} \quad z \longrightarrow 0 \quad
\mbox{on} \quad {\mathbb H},$$

\noindent where $\alpha_n \in {\mathbb R}_{> 0}$ with $\alpha_n
\nearrow \infty$, $P_n \in {\mathbb C} [z]$ monic and $a_n \in
{\mathbb C}$.

\medskip\noindent
Moreover, if $\sphericalangle  D / \pi \in {\mathbb R} \setminus
{\mathbb Q}$, then $P_n = 1$ for all $n \in {\mathbb N}_0$. Note
that $\alpha_0 = \sphericalangle  D/\pi$ and $P_0 = 1$ for any
angle. In particular we see that the Riemann map is not subanalytic
if $\sphericalangle D/\pi \in {\mathbb R} \setminus {\mathbb Q}$. To
use this asymptotic development we want to have a quasianalytic
property; we want to realize these Riemann maps in a class of
functions with an asymptotic development as in ($\dagger$) such that
the functions in this class are determined by the (in general not
convergent) asymptotic expansion. Such quasianalyticity properties
are key tools in generating o-minimal structures (see [21], Van den
Dries-Speissegger [11, 12] and Rolin et al. [30]; see also Badalayan
[1] for quasianalytic classes of this kind).

\noindent Exactly the same kind of asymptotic development occurs at
a transition map of a real analytic vector field on ${\mathbb R}^2$
at a hyperbolic singularity (see Ilyashenko [18]). Poincar\'e return
maps are compositions of finitely many transition maps and are an
important tool to understand qualitatively the trajectories and
orbits of a polynomial or analytic vector field on the plane.
Following Dulac's approach (see [13]), Ilyashenko uses asymptotic
properties of the Poincar\'e maps to solve Dulac's problem (the weak
form of (the second part) of Hilbert's 16th problem): a polynomial
vector field on the plane has finitely many limit cycles (see
Ilyashenko [19] for an overview of the history of Hilbert~16, part
2). One of the first steps in Ilyashenko's proof is to show that the
transition maps at hyperbolic singularities are in a certain
quasianalytic class. Formulating his result on the Riemann surface
of the logarithm (compare with the introduction of [21] and with
[21, Proposition~2.8]) he proves that the considered transition maps
have a holomorphic extension to certain subsets of the Riemann
surface of the logarithm, so-called standard quadratic domains (see
Section~2 below), such that the asymptotic development holds
there.\\
By doing reflections at analytic arcs infinitely often we are able
to extend the Riemann map from the upper half plane to a simply
connected domain with an analytic corner (at 0) to a standard
quadratic domain such that the asymptotic development holds there.
As a consequence we can show the following

\bigskip\noindent {\bf Theorem B}

\noindent {\it Let $\Omega \subset {\mathbb C}$ be a bounded simply
connected domain which is semianalytic. Let $\Phi \colon {\mathbb H}
\to \Omega$ be a Riemann map such that 0 is mapped to a boundary
point of $\Omega$ with angle different from 0. Then $\Phi$ can be
realized in the quasianalytic class of Ilyashenko described above.}

\medskip\noindent
Transition maps at a hyperbolic singularity exhibit a similar
dichotomy of the asymptotic development as indicated in ($\dagger$),
depending whether the hyperbolic singularity is resonent or
non-resonant, i.e. whether the ratio of the two eigenvalues of the
linear part of the vector field at the given hyperbolic singularity
is rational or irrational, see [13] and [21]. In [21] it is shown
that transition maps at non-resonant hyperbolic singularities are
definable in a common o-minimal structure, denoted by ${\mathbb
R}_{{\mathcal Q}}$. This is obtained by proving that the functions
(restricted to the positive line) in Ilyashenko's quasianalytic
class which have no log-terms in their asymptotic expansion,
generate an o-minimal structure ${\mathbb R}_{{\mathcal Q}}$. Using
this result we can derive Theorem~A from Theorem~B (with the
o-minimal structure ${\mathbb R}_{{\mathcal Q}}$). Theorem~A (and B)
are also generalized to the case of unbounded simply connected and
globally semianalytic proper domains.

\medskip\noindent
One may think of possible generalizations of Theorem~A. In the case
of angle different from 0 which is a rational multiple of $\pi$ the
Riemann map is realized in the quasianalytic class of Ilyashenko by
Theorem~B. But so far there is no proof of o-minimality if the
asymptotic expansion has logarithmic terms. If the angle is 0 there
is no asymptotic development known, although some asymptotic
behaviour is known. The proof of Theorem~B could be generalized to
domains definable in other o-minimal structures (for example the
o-minimal structure ${\mathbb R}_{\an}^*$, see [11]), but we lose
then the dichotomy of the asymptotic development depending whether
the angle divided by $\pi$ is irrational or rational, which is
essential for the formulation of Theorem~A and for the application
of the results in [21].

\medskip\noindent
This paper is organized as follows: first, we present in a
preliminary section the basic facts about (globally) semi- and
subanalytic sets, o-minimal structures, and the Riemann Mapping
Theorem which we use throughout the text. In Section~1 we rigorously
define what we mean by angle of the boundary at a boundary point of
a semianalytic domain and we introduce the concept of domains with
an analytic corner. In Section~2 we introduce the quasianalytic
class established by Ilyashenko and prove Theorem~B in several
steps, the main one given by domains with an analytic corner. In
Section~3 we obtain Theorem~A and give applications.

\bigskip\noindent
{\bf Notation.}

\noindent By ${\mathbb N}$ we denote the set of natural numbers and
by ${\mathbb N}_0$ the set of nonnegative integers. Let $a \in
{\mathbb C}$ and $r > 0$ . We set $B(a, r) := \{ z \in {\mathbb C}
\vert \; \vert z -a \vert < r \}$ and $\overline{B} (a, r) := \{ z
\in {\mathbb C} \vert \; \vert z - a \vert \leq r \}$, where $\vert
\quad \vert$ is the euclidean norm. A {\it domain} is an open,
nonempty and connected set (in a topological space). A domain in
${\mathbb C}$ is called {\it simply connected} if its complement has
no bounded connected components. So a bounded domain in ${\mathbb
C}$ is simply connected iff its complement is connected. Given an
open set $U$ of a Riemann surface we denote with ${\mathcal O} (U)$
the ${\mathbb C}$-algebra of holomorphic functions on $U$ with
values in ${\mathbb C}$. We identify ${\mathbb C}$ with ${\mathbb
R}^2$.

\newpage\noindent {\large \bf Preliminaries}

\bigskip\noindent
a) {\bf Semi- and subanalytic sets}

\medskip\noindent
A subset $A$ of ${\mathbb R}^n$, $n \geq 1$, is called {\it
semianalytic} if the
following holds:\\
for each $x_0 \in {\mathbb R}^n$ there are open neighbourhoods $U,
V$ of $x_0$ with $\overline{U} \subset V$ and there are
real-analytic functions $f_i, g_{i, 1}, \dots, g_{i, k}$ on $V$, $1
\leq i \leq \ell$, such that
$$A \cap U = \bigcup\limits_{1 \leq i \leq \ell} \{ x \in U \; \vert \; f_i
(x) = 0, \; g_{i, 1} (x) > 0, \dots, g_{i, k} (x) > 0 \}.$$

\noindent A subset $B$ of ${\mathbb R}^n$, $n \geq 1$, is called
{\it subanalytic} if the
following holds:\\
for each $x_0 \in {\mathbb R}^n$ there is an open neighbourhood $U$
of $x_0$, some $m \geq n$ and some bounded semianalytic set $A
\subset {\mathbb R}^m$ such that $B \cap U = \pi_n (A)$ where $\pi_n
\colon {\mathbb R}^m \to {\mathbb R}^n$, $(x_1, \dots, x_m) \mapsto
(x_1, \dots, x_n)$, is the projection on the first $n$ coordinates.

\noindent A map is called semianalytic resp. subanalytic if its
graph is a semianalytic resp. subanalytic set. A set is called {\it
globally} semianalytic resp. {\it globally} subanalytic if it is
semianalytic resp. subanalytic after applying the semialgebraic
homeomorphism ${\mathbb R}^n \to ] -1, 1 [^n$, $x_i \mapsto
x_i/\sqrt{1+x_i^2}$, (or equivalently if it is semianalytic resp.
subanalytic in the ambient projective space, see [7] and [10,
pp.505-506]).

\noindent Semi- and subanalytic sets exhibit nice `tame' behaviour
(see for example [2], [6], [24], and [31]). One-dimensional
(globally) subanalytic sets and (globally) subanalytic subsets of
${\mathbb R}^2$ are (globally) semianalytic (see [2, Theorem 6.1]).
A bounded semianalytic function on the positive real line is given
locally at 0 by a convergent Puiseux series $\sum\limits_{n =
0}^{\infty} a_n t^{\frac{n}{d}}$ for some $d \in {\mathbb N}$ (see
for example [7, p.192]).

\bigskip\noindent
b) {\bf O-minimal structures}

\medskip\noindent O-minimal structures are axiomatically defined as follows.\\
For $n \in {\mathbb N}$ let $M_n$ be a set of subsets of ${\mathbb
R}^n$ and let ${\mathcal M} : = ({\mathcal M}_n)_{n \in {\mathbb
N}}$. Then ${\mathcal M}$ is called a structure on ${\mathbb R}$, if
the following axioms hold for each $n \in {\mathbb N}$:

\begin{itemize}
\item[(S1)] ${\mathcal M}_n$ is a boolean algebra of subsets of
${\mathbb R}^n$ with ${\mathbb R}^n \in {\mathcal M}_n$.
\item[(S2)] If $A \in {\mathcal M}_n$, then ${\mathbb R} \times A$
and $A \times {\mathbb R}$ belong to ${\mathcal M}_{n+1}$.
\item[(S3)] $\{ (x_1, \dots, x_n) \in {\mathbb R}^n \; \vert \; x_1
= x_n \} \in {\mathcal M}_n$ for $n > 0$.
\item[(S4)] If $A \in {\mathcal M}_{n+1}$, then $\pi (A) \in
{\mathcal M}_n$, where $\pi \colon {\mathbb R}^{n+1} \to {\mathbb
R}^n$ is the projection on the first $n$ coordinates.
\end{itemize}

\noindent An o-minimal structure on ${\mathbb R}$ is a structure
${\mathcal M}$ on ${\mathbb R}$ with the additional properties\\
(O0) $\{ r \} \in {\mathcal M}_1$ for all $r \in {\mathbb R}$.\\
(O1) $\{ (x, y) \in {\mathbb R}^2 \; \vert \; x < y \} \in {\mathcal
M}_2$.\\
(O2) The sets in ${\mathcal M}$ are exactly the finite unions of
intervals and points.

\medskip\noindent
Here $<$ is the canonical order on ${\mathbb R}$. An o-minimal
structure ${\mathcal M}$ on ${\mathbb R}$ {\it expands the field}
${\mathbb R}$ if the following holds:\\
(R1) $\{ (x, y, z) \in {\mathbb R}^3 \; \vert \; z = x + y \} \in
{\mathcal M}_3$.\\
(R2) $\{ (x, y, z) \in {\mathbb R}^3 \; \vert \; z = x \cdot y \}
\in {\mathcal M}_3$.

\medskip\noindent
A subset of ${\mathbb R}^n$ is called {\it definable} in the
structure ${\mathcal M}$ on ${\mathbb R}$ if it belongs to $M_n$. A
function is {\it definable} in ${\mathcal M}$ if its graph is
definable in ${\mathcal M}$. Axiom (O2) implies that a subset of
${\mathbb R}$, definable in an o-minimal structure on ${\mathbb R}$,
has finitely many components. Axiom (R1) and (R2) imply that
addition and multiplication are definable in an o-minimal structure
expanding the field ${\mathbb R}$.

\noindent Not only definable subsets of ${\mathbb R}$ have finitely
many connected components, much more can be deduced from the axioms
of o-minimality: a definable subset of ${\mathbb R}^n$, $n \in
{\mathbb N}$ arbitrary, has finitely many connected components that
are again definable. If the o-minimal structure expands the field
${\mathbb R}$ then definable sets can be definably triangulated and
have a definable $C^k$-stratification for any $k \in {\mathbb N}$.
Hence for a given $k \in {\mathbb N}$, a definable function is $C^k$
outside a definable set of small dimension. General facts about
o-minimal structures can be found in [8].

\bigskip\noindent
{\bf Examples of o-minimal structures on the field ${\mathbb R}$:}

\noindent
\begin{itemize}
\item[(i)] There is a `smallest' (with respect to inclusion of
the boolean algebras of definable sets in any dimension) o-minimal
structure on the field ${\mathbb R}$ denoted by $({\mathbb R}, +,
\cdot, <)$. The definable sets are exactly the semialgebraic sets,
i.e. a definable subset of ${\mathbb R}^n$ is a finite union of sets
of the form
$$\{ x \in {\mathbb R}^n \; \vert \; f(x) = 0, \; g_1 (x) > 0,
\dots, g_r (x) > 0 \}$$

\noindent with $f, g_1, \dots, g_r \in {\mathbb R} [X_1, \dots,
X_n]$ (see [3] for more details).

\end{itemize}

\medskip
\noindent Given a structure ${\mathcal M}$ and functions $f_j \colon
{\mathbb R}^{n_j} \to {\mathbb R}$ ($j$ is in some index set $J$),
we denote by ${\mathcal M} ((f_j)_{j \in J})$ the `smallest'
structure which contains all sets definable in ${\mathcal M}$ and
the graphs of all functions $f_j$.

\medskip\noindent
\begin{itemize}
\item[(ii)] ${\mathbb R}_{\exp} = ({\mathbb R}, +, \cdot, <, \exp)$
where $\exp$ is the exponential function $\exp \colon {\mathbb R}
\to {\mathbb R}_{> 0}$ (see [33] for more details).
\item[(iii)] ${\mathbb R}_{\an} = ({\mathbb R}, +, \cdot, <, (f))$ where $f$
ranges over all the {\it restricted analytic} functions. A function
$f \colon {\mathbb R}^n \to {\mathbb R}$ is called restricted
analytic if the following holds:
$$f = \left\{ \begin{array}{ll}
g \quad \mbox{on} \quad [-1, 1]^n, & g \quad \mbox{analytic on a neighbourhood} \\
                                   & \mbox{of} \quad [-1, 1]^n, \\
0 \quad \mbox{outside} \quad [-1, 1]^n . &              \\
\end{array} \right.$$

\noindent The sets definable in ${\mathbb R}_{\an}$ are exactly the
globally subanalytic sets and the bounded sets definable in
${\mathbb R}_{\an}$ are exactly the bounded subanalytic sets (see
[7] and [10, p.505] for more details).
\item[(iv)] ${\mathbb R}_{\an}^{{\mathbb R}} = {\mathbb R}_{\an} ((x^r)_{r \in {\mathbb
R}})$ where $x^r$ is given by
$$x^r \colon {\mathbb R} \longrightarrow {\mathbb R}, \; x \longmapsto
\left\{ \begin{array}{ccl} a^r & , & a > 0, \\
0 & , & a \leq 0, \\
\end{array} \right.$$

\noindent (see Miller [25] for more details).
\item[(v)] ${\mathbb R}_{\an, \exp} = {\mathbb R}_{\an} (\exp)$ where exp is the
exponential function $\exp \colon {\mathbb R} \to {\mathbb R}_{> 0}$
(see Van den Dries et al. [9] for more details).

\item[(vi)] ${\mathbb R}^*_{\an}$, the o-minimal structure
in which convergent generalized real power series are definable (see
[11] for more details).
\item[(vii)] ${\mathbb R}_{{\mathcal Q}}$, the o-minimal structure in
which transition maps of real analytic vector fields on the plane at
non-resonant hyperbolic singularities are definable (see [21] for
more details).
\end{itemize}

\bigskip\noindent c) {\bf Riemann Mapping Theorem}

\noindent {\it Let $\Omega \subsetneqq {\mathbb C}$ be a simply
connected domain. Then $\Omega$ can be mapped biholomorphically onto
the unit ball} $B(0, 1) = \{ z \in {\mathbb C} \; \vert \; \vert z
\vert < 1 \}$.

\medskip\noindent
Such a biholomorphic map (from a simply connected domain onto the
unit ball) is not totally unique, but the following holds (see [29,
p.179]):

\bigskip\noindent
{\bf Uniqueness}

\noindent Let $a \in \Omega$. Then there is exactly one
biholomorphic map $\varphi \colon \Omega \to B (0, 1)$ with $\varphi
(a) = 0$ and $\varphi' (a) > 0$.

\bigskip\noindent {\bf Examples}

\begin{itemize}
\item[(i)] The group of holomorphic automorphisms of the unit ball is
given by
$$\mbox{Aut} \; B(0, 1) = \left\{ z \longmapsto \rho
\frac{z-a}{\overline{a} z-1} \; \big\vert \; \vert a \vert < 1, \;
\vert \rho \vert = 1 \right\}.$$
\item[(ii)] Let ${\mathbb H} : = \{ z \in {\mathbb C} \; \vert \;
\mbox{Im}\, z > 0 \}$ be the upper half plane. Then ${\mathbb H}
\stackrel{\cong}{\longrightarrow} B (0, 1)$, $z \longmapsto
\frac{z-i}{z+i}$, is a biholomorphic map onto the unit ball.
\end{itemize}

\medskip\noindent
The maps from the example above are all fractional linear maps,
so-called M\"obius transformations (see Conway [5, III. \S 3]). With
the uniqueness result above we see that the M\"obius transformations
are exactly the biholomorphic maps between sets which are either
balls or open half planes. M\"obius transformations are
semialgebraic (identifying ${\mathbb C}$ with ${\mathbb R}^2$) and
are therefore definable in every o-minimal structure (expanding the
field ${\mathbb R}$). A biholomorphic map from one simply connected
domain onto a second one will be here referred to as a Riemann map.
Having one Riemann map from a simply connected domain onto the unit
ball you get all of them by composing with the M\"obius
transformations, which are
automorphisms of the unit ball.\\
Let $\Omega \subset {\mathbb C}$ be a bounded and semianalytic
domain. Then since the o-minimal structure ${\mathbb R}_{\an}$ has
analytic cell decomposition (see [10, pp.508-509]), the boundary is
real analytic at all but finitely many boundary points. Given a
boundary point $x \in
\partial \Omega$ there is $k \in {\mathbb N}$ such that for all $r >
0$ small enough, $\Omega \cap B(x, r)$ consists of $k$ connected
components (each semianalytic) such that each of them is a Jordan
domain. A Jordan domain is a domain whose boundary is a closed
Jordan curve. A Jordan domain is simply connected. Let $C$ be one of
the connected components of $\Omega \cap B (x, r)$, $r > 0$ small.
If $\Omega$ is simply connected and $\varphi$ a Riemann map onto a
simply connected bounded and semianalytic domain $\Omega'$, then by
the Curve Selection Lemma (see [8, p.94]) and Carath\'eodory's Prime
End Theorem (see [28, Chapter 2 p.18]), $\varphi$ has a continuous
extension to $\overline{C}$ with $\varphi (\overline{C} \cap
\partial \Omega) \subset \partial \Omega'$.

\newpage\noindent {\large \bf 1. Angles and domains with an analytic
corner}

\bigskip\noindent
Let $\Omega$ be a bounded and subanalytic domain in ${\mathbb R}^n$.
Let $x \in \partial \Omega := \overline{\Omega} \setminus \Omega$.
Then the germ of $\Omega$ at $x$ has finitely many connected
components. More precisely we have the following: there is $k \in
{\mathbb N}$ such that for all sufficiently small neighbourhoods $V$
of $x$ the set $\Omega \cap V$ has exactly $k$ components having $x$
as boundary
point.\\
Let $\Omega \subset {\mathbb R}^2$ be a bounded and semianalytic
domain without isolated boundary points. Let $x \in
\partial \Omega$ and let $C$ be a connected component of the germ of
$\Omega$ at $x$. Then the germ of the boundary of $C$ at $x$ is
given by (the germs of) two semianalytic curves. So the interior
angle of $\partial C$ at $x$, denoted by $\sphericalangle_{x} C$, is
well defined; it takes value in $[0, 2 \pi]$. If the germ of
$\Omega$ at $x$ is connected we write $\sphericalangle_{x} \Omega$.

\bigskip\noindent
{\bf Definition 1.1.}

\noindent Let $\Omega \subset {\mathbb R}^2$ be a bounded
semianalytic domain without isolated boundary points.
\begin{itemize}
\item[a)] A point $x \in \partial \Omega$ is a {\it singular boundary point}
if $\partial \Omega$ is not a real analytic manifold at $x$.
\item[b)] We set $Sing (\partial \Omega) := \{ x \in \partial \Omega \;
\vert \; x$ is a singular boundary point of $\partial \Omega \}$.
\item[c)] Let $x \in \partial \Omega$. We set $\sphericalangle (\Omega, x) := \{
\sphericalangle_{x} C \; \vert \; C$ is a component of the germ of
$\Omega$ at $x$ and $x \in \Sing (\partial C) \}$.
\end{itemize}

\bigskip\noindent {\bf Remark 1.2.}

\noindent Let $\Omega \subset {\mathbb R}^2$ be a bounded
semianalytic domain without isolated boundary points.
\begin{itemize}
\item[a)] Then $\Sing (\partial \Omega)$ is finite by analytic cell
decomposition (see [10, pp.508-509]).
\item[b)] Let $x \in \partial \Omega$. Then $\sphericalangle (\Omega, x) = \emptyset$
iff $x \not\in \Sing (\partial C)$ for all components $C$ of the
germ of $\Omega$ at $x$. This is especially the case if $x \not\in
\Sing (\partial \Omega)$.
\end{itemize}

\bigskip\noindent {\bf Example 1.3.}
\begin{itemize}
\item[a)] We consider $\Omega := B(0, 1) \setminus \left[ - \frac{1}{2},
\frac{1}{2} \right]$. Then $\partial \Omega = \partial B (0, 1) \cup
\left[ - \frac{1}{2}, \frac{1}{2} \right]$ and $\Sing (\partial
\Omega) = \left\{ - \frac{1}{2}, \frac{1}{2} \right\}$. Given $x \in
\partial B(0, 1)$ the germ of $\Omega$ at $x$ has one component and
we have $\sphericalangle_{x} \Omega = \pi$. Given $x \in \left] -
\frac{1}{2}, \frac{1}{2} \right[$ the germ of $\Omega$ at $x$ has
two components $C_1, C_2$, and we have $\sphericalangle_{\; x} C_1 =
\sphericalangle_{x}C_2 = \pi$. For $x = \pm \frac{1}{2}$ the germ of
$\Omega$ at $x$ has one component $C_{\pm}$ with
$\sphericalangle_{\pm \frac{1}{2}} C_{\pm} = 2 \pi$. Hence
$\sphericalangle \left( \Omega, \pm \frac{1}{2} \right) = \{ 2 \pi
\}$.
\item[b)] We consider the simply connected domain $\Omega : = B (1, 1) \setminus \overline{B}\!\left(
\frac{1}{2}, \frac{1}{2} \right)$. Then $0 \in \Sing (\partial
\Omega)$ and the germ of $\Omega$ at 0 has two components. We have
$\sphericalangle (\Omega, 0) = \{ 0 \}$.
\end{itemize}

\bigskip\noindent
{\bf Definition 1.4.}

\noindent Let $\Omega \subsetneqq {\mathbb R}^2$ be a (not
necessarily bounded) globally semianalytic domain (or equivalently,
definable in the o-minimal structure ${\mathbb R}_{\an}$) without
isolated boundary points.
\begin{itemize}
\item[a)] By $\partial^{\infty} \Omega$ we denote the boundary of $\Omega$
with respect to the standard topology in ${\mathbb R}^2 \cup \{
\infty \}$.
\item[b)] We set $\Omega' := \frac{1}{\Omega} \cap B (0, 1)$, where
$\frac{1}{\Omega} := \left\{ \frac{1}{z} \; \vert \; z \in \Omega
\setminus \{ 0 \} \right\}$.
\begin{itemize}
\item[(i)] We define $Sing(\partial \Omega)$ as follows:
Let $x \in \partial^{\infty} \Omega$. If $x \ne \infty$ then $x \in
Sing (\partial \Omega)$ iff $x \in \Sing (\partial (\Omega \cap B
(0, \vert x \vert +1)))$. If $x = \infty$ we have $\infty \in Sing
(\partial \Omega)$ iff $0 \in \Sing (\partial \Omega')$.
\item[(ii)] Let $x \in \partial^{\infty} \Omega$ and let $C$ be a component of the germ of
$\Omega$ at $x$. If $x \ne \infty$ then $C$ is also a germ of
$\Omega \cap B (0, \vert x \vert + 1)$ at $x$ and we define
$\sphericalangle_{x} C$ as in Remark 1.2. If $x = \infty$ we set
$\sphericalangle_{\infty} C := \sphericalangle_{0} C'$ with $C' :=
\frac{1}{C} \cap B (0,1)$.
\item[(iii)] Let $x \in \partial^{\infty} \Omega$. If $x \ne \infty$ we set
$\sphericalangle (\Omega, x) := \sphericalangle (\Omega \cap B (0,
\vert x \vert +1))$. If $x = \infty$ we set $\sphericalangle
(\Omega, \infty) := \sphericalangle (\Omega', 0)$.
\end{itemize}
\end{itemize}

\bigskip\noindent
Let $\Omega$ be a semianalytic and simply connected bounded domain.
Let $\Phi \colon {\mathbb H} \to \Omega$ be a Riemann map. By
Carath\'eodory's Prime End Theorem we know that $\Phi$ has a
continuous extension to $\overline{{\mathbb H}}$ (see part c) of the
preliminary section). Let $x := \Phi (0) \in
\partial \Omega$. Then there is a component $C$ of the germ of
$\Omega$ at $x$ such that the germ of ${\mathbb H}$ at $0$ is mapped
conformally to $C$ by $\Phi$. We say that 0 is mapped by $\Phi$ to
$x$ with attached angle $\sphericalangle_{x} C$.

\bigskip\noindent
{\bf Definition 1.5} (compare with [23]).

\noindent We say that a domain $D \subset {\mathbb C}$ with $0 \in
\partial D$ has an analytic corner (at 0) if the boundary of $D$ at
0 is given by two analytic arcs which are regular at 0 and if $D$
has an interior angle greater than 0.
More precisely, the following holds:\\
There are holomorphic functions $\varphi_{(1)}, \varphi_{(2)} \in
{\mathcal O} (B(0,1))$ with $\varphi_{(1)} (0) = \varphi_{(2)} (0) =
0$ and $\varphi'_{(1)} (0) \cdot \varphi'_{(2)} (0) \ne 0$ such that
with $\Gamma_1 := \varphi_{(1)} ([0, 1[)$ and $\Gamma_2 :=
\varphi_{(2)} ([0, 1[)$ the following holds:

\medskip\noindent
a) $\partial D \cap B (0, r) = (\Gamma_1 \cup \Gamma_2) \cap B (0,
r)$ for some $r > 0$.

\noindent b) The interior angle $\sphericalangle D \in [0, 2 \pi ]$
of $\partial D$ at 0 is greater than 0.

\medskip\noindent
Note that possibly $\Gamma_1 = \Gamma_2$ if $\sphericalangle D = 2
\pi$. Otherwise we may assume that $\Gamma_1 \cap \Gamma_2 = \{ 0
\}$.

\bigskip\noindent {\bf Remark 1.6.}

\noindent Let $D \subset {\mathbb C}$ be a simply connected domain
with an analytic corner at $0 \in \partial D$. Let $\Phi \colon
{\mathbb H} \to D$ be a Riemann map. Again by Carath\'eodory's Prime
End Theorem (see part c) of the preliminary section) we get the
following: the inverse $\Phi^{-1}$ has a continuous extension to the
boundary $\partial D$ in a neighbourhood of $0 \in
\partial D$ with $\Phi^{-1} (0) = \colon x \in
\partial {\mathbb H}$ and $\Phi$ has a continuous extension to the boundary
$\partial {\mathbb H}$ in a neighbourhood of $x \in \partial
{\mathbb H}$ with $\Phi (x) = 0$. Choosing an automorphism $\tau$ of
${\mathbb H}$ (a certain M\"obius transformation) which maps 0 to
$x$ and replacing $\Phi$ by $\Phi \circ \tau$ we can assume that
$\Phi (0) = 0$. Renaming the indices $i = 1, 2$ from Definition 1.5
we can assume that the germ of ${\mathbb R}_{\geq 0}$ at 0 is mapped
by $\Phi$ to the germ of $\Gamma_1$ at 0 and that the germ of
${\mathbb R}_{\leq 0}$ at 0 is mapped by $\Phi$ to the germ of
$\Gamma_2$ at 0. We say that the positive direction is mapped to
$\Gamma_1$ and the negative direction is mapped to $\Gamma_2$.

\vspace{1.5cm}\noindent {\large \bf 2. Riemann maps at semianalytic
domains and a quasianalytic class of Il\-ya\-shen\-ko}

\medskip\noindent
Lehman showed in [23] that the Riemann map has an asymptotic
development at an analytic corner (see also [28, p.58]; Wigley [32]
showed the existence of asymptotic development for more general
corners). We introduce the series which occur as asymptotic
expansions.

\bigskip\noindent
{\bf Definition 2.1} (compare with [21, Section 1]).

\noindent Let $z$ be an indeterminate. A generalized log-power
series in $z$ is a formal expression $g(z) = \sum\limits_{\alpha \in
{\mathbb R}_{\geq 0}} a_{\alpha} P_{\alpha} (\log z) z^{\alpha}$
with $a_{\alpha} \in {\mathbb C}$ and $P_{\alpha} \in {\mathbb C}
[z] \setminus \{ 0 \}$ monic with $P_0 = 1$ such that the support of
$g$, defined as $\supp (g) : = \{ \alpha \in {\mathbb R}_{\geq 0} \;
\vert \; a_{\alpha} \ne 0 \}$, fulfils the following condition: for
all $R > 0$ the set $\supp (g) \cap [0, R]$ is finite. We write
${\mathbb C} [[z^*]]_{\log}^{\omega}$ for the set of generalized
log-power series. For $g \in {\mathbb C} [[z^*]]_{\log}^{\omega}$ we
set $\nu (g) := \min \supp (g)$. By ${\mathbb C} [[z^*]]^{\omega}$
we denote the subset of ${\mathbb C} [[z^*]]_{\log}^{\omega}$
consisting of all $g \in {\mathbb C} [[z^*]]_{\log}^{\omega}$ with
$P_{\alpha} = 1$ for all $\alpha \in {\mathbb R}_{\geq 0}$. By
${\mathbb C} [[z^*]]_{\log}^{\omega, \fin}$ (resp. ${\mathbb C}
[[z^*]]^{\omega, \fin}$) we denote the set of all $ g \in {\mathbb
C} [[z^*]]^{\omega}_{\log}$ (resp. ${\mathbb C} [[z^*]]^{\omega}$)
with finite support.

\bigskip\noindent
{\bf Convention.} From now on we omit the superscript $\omega$.

\bigskip\noindent
{\bf Remark 2.2.}
\begin{itemize}
\item[a)] The set ${\mathbb C}
[[z^*]]_{\log}$ is in a natural way a ${\mathbb C}$-algebra with
${\mathbb C} [[z^*]]$ as subalgebra.
\item[b)] Interpreting $\log$ as the principal branch of the
logarithm, i.e. as $\log \colon {\mathbb C} \setminus {\mathbb
R}_{\leq 0} \to {\mathbb C}$, $z \mapsto \log \vert z \vert + i \Arg
z$, (with $\Arg (z) \in ] - \pi, \pi [$ the standard argument) and
$z^{\alpha}$ as the power function $z^{\alpha} \colon {\mathbb C}
\setminus {\mathbb R}_{\leq 0} \to {\mathbb C}$, $z \mapsto \exp
(\alpha \log z)$, we get that $g \in {\mathcal O} ({\mathbb C}
\setminus {\mathbb R}_{\leq 0})$ if $g \in {\mathbb C}
[[z^*]]_{\log}^{\fin}$.
\end{itemize}

\bigskip\noindent {\bf Definition 2.3.}

\noindent Let $f \in {\mathcal O} ({\mathbb H})$ and let $g = \Sigma
\; a_{\alpha} P_{\alpha} (\log z) z^{\alpha} \in {\mathbb C}
[[z^*]]_{\log}$. We say that $f$ has asymptotic expansion $g$ on
${\mathbb H}$ and write $f \sim_{{\mathbb H}} g$, if for each $R >
0$
$$f(z) - \sum\limits_{\alpha \leq R} a_{\alpha} P_{\alpha} (\log z)
z^{\alpha} = o (\vert z \vert^R) \quad \mbox{as} \quad \vert z \vert
\longrightarrow 0 \quad \mbox{on}\quad {\mathbb H}.$$

\noindent Note that $g$ is unique.

\bigskip\noindent {\bf Fact 2.4} (see Lehman [23, Theorem 1]){\bf .}

\noindent {\it Let $D$ be a simply connected domain with an analytic
corner and let} $\Phi \colon {\mathbb H} \to D$ {\it be a Riemann
map with} $\Phi (0) = 0$. {\it Then there is} $g \in {\mathbb C}
[[z^*]]_{\log}$ {\it such that} $f \sim_{{\mathbb H}} g$. {\it If}
$\sphericalangle  D / \pi \in {\mathbb R} \setminus {\mathbb Q}$
{\it then} $g \in {\mathbb C} [[z^*]]$.

\medskip\noindent
In the above situation we have with $\alpha := \sphericalangle
D/\pi$ that $\supp (g) \subset {\mathbb N}_0 + {\mathbb N} \alpha$,
that $\nu(g) = \alpha$ and that $P_{\nu (g)} = 1$. Moreover, by
Remark 1.6 the Riemann map has a continuous extension to the
boundary of ${\mathbb H}$ near 0. The estimates of Definition~2.3
then hold on $\overline{{\mathbb H}} \cap B(0, r)$ for some $r
> 0$ (with the continuous extension of the logarithm and the power
functions to $\overline{{\mathbb H}}$ respectively), see [23,
Section 4].

\medskip\noindent
We want to realize the Riemann map in the quasianalytic class used
by Ilyashenko in his work on Hilbert~16. For this we have to
consider holomorphic functions on the Riemann surface of the
logarithm (compare with [21, Section 2]):

\bigskip\noindent {\bf Definition 2.5.}

\noindent We define the Riemann surface of the logarithm ${\mathbf
L}$ in polar coordinates by ${\mathbf L} := {\mathbb R}_{> 0} \times
{\mathbb R}$. Then ${\mathbf L}$ is a Riemann surface with the
isomorphic holomorphic projection map $\log \colon {\mathbf L} \to
{\mathbb C}$, $(r, \varphi) \mapsto \log r + i \varphi$. For $z =
(r, \varphi) \in {\mathbf L}$ we define the absolute value by $\vert
z \vert := r$ and the argument by $\arg z := \varphi$. For $r
> 0$ we set $B_{{\mathbf L}}(r) := \{ z \in {\mathbf L} \; \vert \; \vert z
\vert < r \}$. We identify ${\mathbb C} \setminus {\mathbb R}_{\leq
0}$ with ${\mathbb R}_{> 0} \times ] - \pi, \pi [ \subset {\mathcal
L}$ via polar coordinates. Let $\alpha \geq 0$. We define the power
function $z^{\alpha}$ as $z^{\alpha} \colon {\mathbf L} \to {\mathbb
C}$, $z = (r, \varphi) \mapsto \exp (\alpha \log z)$. For each real
$\rho > 0$, the map ${\mathbf p} \colon {\mathbf L} \to {\mathbf L}$
is defined by ${\mathbf p}^{\rho} (r, \varphi) := (r^{\rho}, \rho
\varphi)$, and the map ${\mathbf m} \colon {\mathbf L}^2 \to
{\mathbf L}$ is defined by ${\mathbf m} ((r_1, \varphi_1), (r_2,
\varphi_2)) := (r_1 r_2, \varphi_1 + \varphi_2)$ (see [21, Section
4]).

\newpage\noindent {\bf Definition 2.6.}

\noindent A domain $W \subset {\mathbf L}$ of the Riemann surface of
the logarithm is a standard quadratic domain if there are constants
$c, C > 0$ such that
$$W = \left\{ (r, \varphi) \in {\mathbf L} \; \big\vert \; 0 < r < c \; \exp (-
C \sqrt{\vert \varphi \vert}) \right\}.$$

\noindent A domain is called a quadratic domain if it contains a
standard quadratic domain.

\bigskip\noindent
{\bf Definition 2.7.}

\noindent Let $U \subset {\mathbf L}$ be a quadratic domain, let $f
\in {\mathcal O} (U)$ and let $g = \sum\limits_{\alpha \geq 0}
a_{\alpha} P_{\alpha} (\log z) z^{\alpha} \in {\mathbb C}
[[z^*]]_{\log}$. We say that $f$ has asymptotic expansion $g$ on $U$
and write $f \sim_U g$, if for each $R > 0$ there is a quadratic
domain $U_R \subset U$ such that
$$f(z) - \sum\limits_{\alpha \leq R} a_{\alpha} P_{\alpha} (\log z)
z^{\alpha} = o (\vert z \vert^R) \quad \mbox{as} \quad \vert z \vert
\longrightarrow 0 \quad \mbox{on} \quad U_R.$$

\noindent We write $T f := g$. By ${\mathcal Q}^{\log} (U)$ we
denote the set of all $f \in {\mathcal O} (U)$ with an asymptotic
expansion. By ${\mathcal Q}(U)$ we denote the subset of all $f \in
{\mathcal Q}^{\log} (U)$ with $Tf \in {\mathbb C} [[z^*]]$.

\bigskip\noindent
{\bf Remark 2.8.}

\noindent
\begin{itemize}
\item[a)] Logarithm and power functions on
${\mathbf L}$ extend logarithm and power functions on ${\mathbb C}
\setminus {\mathbb R}_{\leq 0}$ (compare with Remark~2.2). Given $g
\in {\mathbb C} [[z^*]]^{\fin}_{\log}$ we get that $g \in {\mathcal
O} ({\mathbf L})$. Especially in the situation of Definition~2.7 we
get for $R > 0$ that $\sum\limits_{\alpha \leq R} a_{\alpha}
P_{\alpha} (\log z) z^{\alpha} \in {\mathcal O} ({\mathbf L})$.
\item[b)] If $f \in {\mathcal Q}^{\log}(U)$ for some quadratic domain $U$ then there
is exactly one $g \in {\mathbb C} [[z^*]]_{\log}$ with $f \sim_U g$;
i.e. $Tf$ is well defined.
\end{itemize}

\bigskip\noindent
{\bf Definition 2.9.}

\noindent We define an equivalence relation $\equiv$ on
$\bigcup\limits_{U \subset {\mathbf L} \; \mbox{\tiny quadr.}}
{\mathcal Q}^{\log} (U)$ as follows: $f_1 \equiv f_2$ if and only if
there is a quadratic domain $V \subset {\mathbf L}$ such that $f_1
\vert_V = f_2 \vert_V$. We let ${\mathcal Q}^{\log}$ be the set of
all $\equiv$-equivalence classes. In the same way we obtain the
class ${\mathcal Q}$. Note that ${\mathcal Q} = {\mathcal Q}_1^1$ in
the notation of [21, Definition~5.1, Remarks~5.2 \& Definition 5.4].

\bigskip\noindent {\bf Remark 2.10.}
\begin{itemize}
\item[a)] We will not distinguish between $f \in \bigcup\limits_{U \subset
{\mathbf L} \; \mbox{\tiny quadr.}} {\mathcal Q}^{\log} (U)$ and its
equivalence class in ${\mathcal Q}^{\log}$, which we also denote by
$f$. Thus ${\mathcal Q}^{\log} (U) \subset {\mathcal Q}^{\log}$
given a quadratic domain $U \subset {\mathbf L}$.
\item[b)] In the same way we define ${\mathcal Q} \subset {\mathcal Q}^{\log}$.
We have ${\mathcal Q}(U) \subset {\mathcal Q}$ for $U \subset
{\mathbf L}$ a quadratic domain.
\item[c)] Given a quadratic domain $U \subset {\mathbf L}$ the set ${\mathcal Q}^{\log}
(U)$ is a ${\mathbb C}$-algebra with ${\mathcal Q}(U)$ as a
subalgebra. Also, ${\mathcal Q}^{\log}$ is an algebra with
${\mathcal Q}$ as a subalgebra.
\item[d)] Given a quadratic domain $U \subset {\mathbf L}$ the well defined
maps $T \colon {\mathcal Q}^{\log} (U) \to {\mathbb C}
[[z^*]]_{\log}$, $f \mapsto Tf$, and $T \colon {\mathcal Q}(U) \to
{\mathbb C} [[z^*]]$, $f \mapsto Tf$, are homomorphisms of ${\mathbb
C}$-algebras. Also the induced maps $T \colon {\mathcal Q}^{\log}
\to {\mathbb C} [[z^*]]_{\log}$, $f \mapsto Tf$, and $T \colon
{\mathcal Q} \to {\mathbb C} [[z^*]]$, $f \mapsto Tf$, are
homomorphisms of ${\mathbb C}$-algebras.
\end{itemize}

\newpage\noindent {\bf Proposition 2.11.}

\noindent Let $U \subset {\mathbf L}$ be a quadratic domain. The
homomorphism $T \colon {\mathcal Q}^{\log} (U) \to {\mathbb C}
[[z^*]]_{\log}$ is injective. Therefore, the homomorphism $T \colon
{\mathcal Q}^{\log} \to {\mathbb C} [[z^*]]_{\log}$ is injective.

\bigskip\noindent
Proof:

\noindent See Ilyashenko [18, Theorem 2 p.23] and [21, Proposition
2.8]. \hfill $\square$

\bigskip\noindent
Proposition 2.11 contains the necessary quasianalyticity we need for
our o-minimality result. We realize now the Riemann map from the
upper half plane to a bounded semianalytic domain with attached
angle greater than 0 in this quasianalytic class. Theorem~B gets the
following precise form:

\bigskip\noindent {\bf Theorem 2.12.}

\noindent {\it Let $\Omega \subset {\mathbb C}$ be a semianalytic
and simply connected bounded domain. Let $\Phi \colon {\mathbb H}
\to \Omega$ be a Riemann map such that 0 is mapped to a boundary
point with attached angle $\sphericalangle$ greater than 0. Then
there is a quadratic domain $U \subset {\mathbf L}$ such that $\Phi
\in {\mathcal Q}^{\log} (U)$. If $\sphericalangle \in \pi ({\mathbb
R} \setminus {\mathbb Q})$ then $\Phi \in {\mathcal Q}(U)$.}

\medskip\noindent
We will prove this theorem in several steps, considering first
domains with an analytic corner. We obtain the extension of the
Riemann map to a quadratic domain of the Riemann surface of the
logarithm by performing reflections at analytic arcs infinitely
often. These reflections are obtained by iteration, inversion and
conjugation of certain holomorphic functions. To get the desired
properties of Definition~2.7 we have to control very carefully this
discrete dynamical system. Therefore we use the theory of univalent
functions, especially Koebe's $\frac{1}{4}$-Theorem and the Growth
Theorem (see for example Duren [14, Chapter 2]). To motivate the
technical statements of the upcoming proofs we give the following
example for the Schwarz reflection principle at analytic arcs (which
reduces to the Schwarz reflection principle at the real line, see
[5, IX. 1.1]):

\bigskip\noindent {\bf Example 2.13.}

\noindent Let $r > 0$ and let $V := {\mathbb H} \cap B (0, r)$. Let
$f \in {\mathcal O} (V)$ have the following property: $f$ has a
continuous extension to $[0, r[$ and there is an injective
holomorphic function $\varphi \colon B (0, 1) \to {\mathbb C}$ such
that $f([0, r[) \subset \Gamma := \varphi ([0, 1 [)$. Then there is
some $0 < r' \leq r$ such that $f$ has a holomorphic extension to
$B(0, r') \setminus {\mathbb R}_{\leq 0}$, given by
$$
f \colon B (0, r') \setminus {\mathbb R}_{\leq 0} \longrightarrow
{\mathbb C}, \; z \longmapsto \left\{ \begin{array}{ccl}
f(z) & & \mbox{Im}\, z \geq 0, \\
       & \mbox{if} &              \\
\varphi \overline{(\varphi^{-1} f(\overline{z}))} & & \mbox{Im}\, z
< 0.
\end{array} \right.
$$

\bigskip\noindent
{\bf Definition 2.14.}

\noindent Let $k \in {\mathbb N}_0$. We define $T_k := \{ (r,
\varphi) \in {\mathbf L} \; \vert \; 0 \leq \varphi \leq 2^k \pi
\}$. For $k \in {\mathbb N}$ we set $T'_k := \{ (r, \varphi) \in
{\mathbf L} \; \vert \; 2^{k-1} \pi \leq \varphi \leq 2^k \pi \}$.
Given $k \in {\mathbb N}_0$ we define the reflections $\tau_k \colon
T'_{k+1} \to T_k$, $(r, \varphi) \mapsto (r, - \varphi + 2^{k+1}
\pi)$. Note that $T_0 = \overline{{\mathbb H}} \setminus \{ 0 \}$,
$T_{k+1} = T_k \cup T'_{k+1}$ and $\tau_k (r, 2^k \pi) = (r, 2^k
\pi)$, $\tau_k (r, 2^{k+1} \pi) = (r, 0)$.

\bigskip\noindent {\bf Lemma 2.15.}

\noindent a) {\it Let $k \geq 0$. Then $\overline{\log \circ \tau_k}
= \log - i 2^{k+1}
\pi$.}\\
b) {\it Let $\alpha > 0$ and $k \geq 0$. Then there is some $a \in
{\mathbb C}^*$ with $\vert a \vert = 1$ such that
$\overline{z^{\alpha} \circ \tau_k} = a z^{\alpha}$.}

\bigskip\noindent Proof:

\noindent a) Let $z = (r, \varphi) \in T'_{k+1}$, $k \geq 0$. Then
by Definition 2.5
$$\begin{array}{rcl}
\overline{(\log \circ \tau_k) (r, \varphi)} & = & \overline{\log (r,
- \varphi
+ 2^{k+1} \pi)} \\
& = & \log (r, \varphi) - i 2^{k+1} \pi. \\
\end{array}$$

\noindent b) Let $z = (r, \varphi) \in T'_{k+1}$, $k \geq 0$. Then
by Definition 2.5
$$\begin{array}{rcl}
\overline{(z^{\alpha} \circ \tau_k) (r, \varphi)} & = &
\overline{\exp (\alpha (\log r + i (-\varphi
+ 2^{k+1} \pi))} \\
& = & \exp (- i \alpha 2^{k+1} \pi) z^{\alpha} (r, \varphi). \\
\end{array}$$

\noindent \hfill $\square$

\bigskip\noindent
{\bf Lemma 2.16.}

\noindent {\it Let $r > 0$ and let $f \colon B (0, r) \to {\mathbb
C}$ be holomorphic and
injective with $f(0) = 0$ and $\vert f' (0) \vert = 1$.\\
Then the following holds:}\\
a) $f(B(0, r)) \supset B\!\left( 0, \frac{r}{4} \right)$,\\
b) $\vert f(z) \vert \leq 4 \vert z \vert$ {\it for} $z \in
B\!\left( 0, \frac{r}{2} \right)$.

\bigskip\noindent
Proof:

\noindent a) is a consequence of Koebe's $\frac{1}{4}$-Theorem (see
[14, Theorem 2.3]) and b) is a consequence of the Growth Theorem
(see [14, Theorem 2.6]). \hfill $\square$

\bigskip\noindent
{\bf Theorem 2.17.}

\noindent {\it Let $D$ be a simply connected domain with an analytic
corner. Let $\Phi \colon {\mathbb H} \to D$ be a Riemann map with
$\Phi (0) = 0$. Then there is a quadratic domain $U$ such that
$\Phi$ has a holomorphic extension to $U$.}

\bigskip
\noindent Proof:

\noindent Let $\varphi_{(1)}, \varphi_{(2)}$ and $\Gamma_1$,
$\Gamma_2$ be as in Definition 1.5. We may assume that the positive
direction is mapped to $\Gamma_1$ and the negative direction to
$\Gamma_2$. Replacing $\varphi_{(i)} (z)$ by $\varphi_{(i)} \left(
\frac{z}{\vert \varphi'_{(i)} (0) \vert} \right)$ we get the
following.

\noindent There is some $0 < \overline{r} \leq 1$ such that
$\varphi_{(i)} \colon B (0, \overline{r}) \to {\mathbb C}$ is
injective, $\vert \varphi'_{(i)} (0) \vert = 1$ and $\partial D \cap
B\!\left( 0, \frac{\overline{r}}{4} \right) = (\Gamma_1 \cup
\Gamma_2) \cap B\!\left( 0, \frac{\overline{r}}{4} \right)$.\\
Given $k \geq 0$ we recursively define positive real numbers $r_k$
and injective holomorphic functions $\varphi_{k} \in {\mathcal O}
(B(0, r_k))$, $\chi_{k} \in {\mathcal O} (B\!\left(0, \frac{r_k}{8}
\right))$ with $\varphi_{k} (0) = \chi_{k} (0) = 0$ and $\vert
\varphi'_{k} (0) \vert = \vert \chi'_{k} (0) \vert = 1$ as follows:

\medskip\noindent
$k = 0$: We choose $0 < r_0 \leq \overline{r}$. We will specify
$r_0$ below. We take $\varphi_{0} := \varphi_{(2)} \in {\mathcal O}
(B(0, r_0))$ and
$$\chi_{0} \colon B\!\left( 0, \frac{r_0}{8} \right) \longrightarrow
{\mathbb C}, \; z \longmapsto \overline{\varphi_{0}
(\overline{\varphi_{0}^{-1} (z)})}.$$

\noindent Note that $\chi_{0}$ is well defined: by Lemma 2.16 a) we
know that $\varphi_{0}^{-1} \in {\mathcal O} (B\!\left(0,
\frac{r_0}{4} \right))$. By Lemma 2.16 b) we see that $\vert
\varphi_{0}^{-1} (z) \vert < r_0$ for $\vert z \vert <
\frac{r_0}{8}$.

\medskip\noindent
$k \to k+1$: We take $r_{k+1} := \frac{r_k}{32}$ and $\varphi_{k+1}
\colon B (0, r_{k+1}) \to {\mathbb C}$, $z \mapsto
\overline{\chi_{k} (\varphi_{(1)} (\overline{z}))}$. We set
$$\chi_{k+1} \colon B\!\left(
0, \frac{r_{k+1}}{8} \right) \to {\mathbb C}, z \mapsto
\overline{\varphi_{k+1} (\overline{\varphi_{k+1}^{-1} (z)})}.$$

\noindent Note that $\varphi_{k+1}$ is well defined by Lemma 2.16 b)
and that $\chi_{k+1}$ is well defined by the same argument as in the
case $k = 0$.

\medskip
\noindent By Fact 2.4 and the subsequent remark we find some
$\overline{t}
> 0$ and some $E > 1$ such that $\Phi \in {\mathcal O} (\overline{\mathbb H} \cap B
(0, \overline{t}))$ and $\vert \Phi (z) \vert \leq E \vert z
\vert^{\alpha}$ for $z \in \overline{{\mathbb H}} \cap B (0,
\overline{t})$ where $\alpha := \sphericalangle D / \pi$.

\medskip\noindent
For $k \geq 0$ we recursively define positive real numbers $t_k$ and
holomorphic functions $\Phi_{k} \in {\mathcal O}
(\stackrel{\circ}{T}_k \cap B_{{\mathbf L}}(t_k)) \cap C^0 (T_k \cap
B_{{\mathbf L}}(t_k))$ with $\Phi_k (0) = 0$ and $\Phi_k (T_k \cap
B_{{\mathbf L}} (t_k)) \subset B (0, \frac{r_k}{16})$ as follows:

\medskip\noindent
$k = 0$: We choose $0 < t_0 \leq \overline{t}$ such that $\Phi
(\overline{{\mathbb H}} \cap B (0, t_0)) \subset B (0,
\frac{r_0}{16})$. We set $\Phi_{0} := \Phi$.

\medskip\noindent
$k \to k+1$: We define
$$
\Phi_{k+1} \colon T_{k+1} \cap B_{{\mathbf L}} (t_{k})
\longrightarrow {\mathbb C}, \; z \longmapsto \left\{
\begin{array}{lcl}
\Phi_{k} (z) & & z \in T_k, \\
           & \mbox{if} &             \\
\overline{\chi_{k} (\Phi_{k} (\tau_k (z)))} & & z \in T'_{k+1}. \\
\end{array} \right.
$$

\noindent Note that $\Phi_{k+1}$ is well defined by the construction
of $\chi_k$ and the choice of $t_k$. Moreover, $\Phi_{k+1}$ is a
holomorphic extension of $\Phi_k$. We choose $0 < t_{k+1} \leq t_k$
such that $\Phi_{k+1} (T_{k+1} \cap B_{{\mathbf L}} (t_{k+1}))
\subset B (0, \frac{r_{k+1}}{16})$. Note that this is possible since
$\lim\limits_{\vert z \vert \to 0} \; \Phi_{k+1} (z) = 0$.

\medskip
\noindent So far the only condition imposed on $r_0$ is that $0 <
r_0 < \overline{r}$. We choose now $r_0$ so small such that $\left(
\frac{r_0}{E \cdot 16} \right)^{\frac{1}{\alpha}} \leq
\overline{t}$. For $k \geq 0$ we set $E_k := E 4^k$. By induction on
$k \geq 0$ we show that we can choose $t_k : = \left(
\frac{r_k}{E_{k} \cdot 16} \right)^{\frac{1}{\alpha}}$ and that
$\vert \Phi_{k} (z) \vert \leq E_k \vert z \vert^{\alpha}$ for $z
\in T_k \cap B_{{\mathbf L}} (t_k)$.

\medskip\noindent
$k = 0$: By the choice of $r_0$ we have $0 < t_0 \leq \overline{t}$.
Moreover, $\vert \Phi_0 (z) \vert \leq E_0 \vert z \vert^{\alpha}$
for $z \in \overline{{\mathbb H}} \cap B_{{\mathbf L}}
(\overline{t_0})$ by the setting. By the definition of $t_0$ and
$E_0$ we see that $\vert \Phi_0 (z) \vert \leq \frac{r_0}{16}$ for
$z \in \overline{{\mathbb H}} \cap B_{{\mathbf L}} (t_0)$.

\medskip\noindent
$k \to k+1$: We obtain, by applying Lemma 2.17 b) to $\chi_{k} \in
{\mathcal O} \left( B\!\left( 0, \frac{r_k}{8} \right) \right)$ and
by the fact that from the inductive hypothesis we have $\Phi_k (T_k
\cap B_{{\mathbf L}} (t_k)) \subset B \left( 0, \frac{r_k}{16}
\right)$, that for $z \in T'_{k+1} \cap B_{{\mathbf L}} (t_k)$.

$$
\vert \Phi_{k+1} (z) \vert  =  \vert \chi_{k} (\Phi_{k} (\tau_k
(z))) \vert  \leq  4 \; \vert \Phi_{k} (\tau_k (z)) \vert \leq 4 E_k
\, \vert z
\vert^{\alpha} = E_{k+1} \vert z \vert^{\alpha}.\\
$$

\medskip\noindent
By the definition of $t_{k+1}$ and $E_{k+1}$ we obtain that $\vert
\Phi_{k+1} (z) \vert \leq \frac{r_{k+1}}{16}$ for $z \in T'_{k+1}
\cap B_{{\mathbf L}} (t_{k+1})$. Since $\Phi_{k+1}$ coincides with
$\Phi_k$ on $T_k$ the claim follows from the inductive hypothesis.

\medskip\noindent
By construction $\Phi_{k+1}$ extends $\Phi_k$ holomorphically for
all $k \geq 0$. Hence $\Phi$ has a holomorphic extension to
$\bigcup\limits_{k \geq 0} T_k \cap B_{{\mathbf L}} (t_k)$. By the
definition of $r_k$ and $E_k$ we obtain some $K > 1$ such that $t_k
\geq K^{-k}$ for all $k \geq 1$. For $\varphi > 0$ let $k(\varphi)
\in {\mathbb N}$ such that $2^{k (\varphi) -1} \pi \leq \varphi \leq
2^{k (\varphi)} \pi$, i.e. ${\mathbb R}_{> 0} \times \{ \varphi \}
\subset T'_{k (\varphi)}$. Then there is some $C > 0$ such that $k
(\varphi) \leq C \log \varphi$ for all $\varphi > 0$ large enough.
Enlarging $K > 1$, if necessary, we get that $\{ (r, \varphi) \in
{\mathbf L} \; \vert \; \varphi > 0$ and $r \leq K^{- \log \varphi}
\} \subset \bigcup\limits_{k \geq 0} T_k \cap B_{{\mathbf L}}
(t_k)$. Repeating the reflection process in the negative direction
we see that $\Phi$ has a holomorphic extension to some quadratic
domain $U$ since $\log \varphi \leq \sqrt{\varphi}$ for $\varphi >
0$. \hfill $\square$

\newpage\noindent {\bf Theorem 2.18.}

\noindent {\it Let $D$ be a simply connected domain with an analytic
corner. Let $\Phi \colon {\mathbb H} \to D$ be a Riemann map with
$\Phi (0) = 0$. Then there is a quadratic domain $U \subset {\mathbf
L}$ such that $\Phi \in {\mathcal Q}^{\log} (U)$. If
$\sphericalangle D/\pi \in {\mathbb R} \setminus {\mathbb Q}$ then
$\Phi \in {\mathcal Q} (U)$.}

\bigskip\noindent Proof:

\noindent We use the results and the notation of the previous proof.
We define $s_k := \min \{ t_k, E_k^{-\frac{2}{\alpha}} \}$. Then
$\Phi_{k} \in {\mathcal O} (\stackrel{\circ}{T}_k \cap B_{{\mathbf
L}} (s_k)) \cap C^0 (T_k \cap B_{{\mathbf L}} (s_k))$ and for $z \in
T_k \cap B_{{\mathbf L}}(s_k)$ we have
$$\vert \Phi_{k} (z) \vert \leq E_k \vert z \vert^{\alpha} \leq E_k \vert
s_k^{\frac{\alpha}{2}} \vert \; \vert z^{\frac{\alpha}{2}} \vert
\leq \vert z \vert^{\frac{\alpha}{2}}. \leqno(1)$$

\noindent Let $g = \sum\limits_{\gamma > 0} a_{\gamma} P_{\gamma}
(\log z) z^{\gamma} \in {\mathbb C} [[z^*]]_{\log}$ with $\Phi
\sim_{{\mathbb H}} g$ (see Fact 2.4). Given $R > 0$ we show that
$\Phi (z) - \sum\limits_{\gamma \leq R} a_{\gamma} P_{\gamma} (\log
z) z^{\gamma} = o (\vert z \vert^R)$ as $\vert z \vert \to 0$ on
some quadratic domain $U_R \subset U$ which completes the proof. We
set
$$\varepsilon_R \colon \overline{{\mathbb H}} \cap B (0, s_0) \to
{\mathbb C}, \; z \mapsto \Phi (z) - \sum\limits_{\gamma \leq R}
a_{\gamma} P_{\gamma} (\log z) z^{\gamma}.$$

\noindent We choose a constant $S_R$ with $R < S_R < \min \{ \gamma
\in \supp (g) \; \vert \; \gamma
> R \}$. Then by Fact 2.4 and the subsequent remark there is some $C_R > 1$ such
that $\vert \varepsilon_R (z) \vert \leq C_R \vert z \vert^{S_R}$ on
$\overline{{\mathbb H}} \cap B (0, s_0)$. Shrinking $S_R$ and $s_0$
(by shrinking the above $r_0$) we can assume that
$$\vert \varepsilon_R (z) \vert \leq \vert z \vert^{S_R} \quad
\mbox{on} \quad \overline{{\mathbb H}} \cap B (0, s_0). \leqno(2)$$

\noindent For $k \geq 0$ we set
$$\varepsilon_{R, k} \colon T_k \cap B_{{\mathbf L}}
(s_k) \to {\mathbb C}, z \mapsto \Phi_{k} (z) - \sum\limits_{\gamma
\leq R} a_{\gamma} P_{\gamma} (\log z) z^{\gamma}.$$

\noindent We fix $R > 0$ and omit this subscript when it is clear
from the context.

\medskip\noindent
Let $\sum\limits_{\ell = 1}^{\infty} c_{k, \ell} z^{\ell}$ be the
power series expansion of $\chi_{k}$ on $B \left( 0, \frac{r_k}{8}
\right)$, $k \geq 0$. Applying Lemma 2.16~b) and using Cauchy's
estimate (see for example [5, IV. 2.14]) we obtain that
$$\vert c_{k, \ell} \vert \leq 4 \left( \frac{16}{r_k} \right)^{\ell-1}
\quad \mbox{for} \quad k \geq 0 \quad \mbox{and} \quad \ell \geq 1.
\leqno(3)$$

\noindent Let $m \in {\mathbb N}$ with $m \alpha/2
> R$ (where $\alpha = \sphericalangle D/\pi$). We define
$$h_{k} \colon B\!\left( 0, \frac{r_k}{8} \right) \longrightarrow {\mathbb
C}, \; z \mapsto \chi_{k} (z) - \sum\limits_{\ell = 1}^m c_{k, \ell}
z^{\ell}.$$

\noindent {\bf Claim 1:} If $z \in B\!\left( 0, \frac{r_k}{16}
\right)$ then
$$\vert h_{k} (z) \vert \leq 4 (m+1) \left( \frac{16}{r_k} \right)^m
\vert z \vert^{m+1}. \leqno(4)$$

\bigskip\noindent
Proof of Claim 1:

\noindent By (3) we get $\big\vert \sum\limits_{\ell = 1}^m c_{k,
\ell} z^{\ell-1} \big\vert \leq 4 m$ on $B\!\left( 0, \frac{r_k}{16}
\right)$. Applying the maximum principle to $\sum\limits_{\ell =
1}^m c_{k, \ell} z^{\ell-1}$ on $B \left( 0, \frac{r_k}{16} \right)$
{we see that $\big\vert \sum\limits_{\ell = 1}^m c_{k, \ell}
z^{\ell} \big\vert \leq 4 m \vert z \vert$ on $B (0,
\frac{r_k}{16})$ and as a consequence we obtain, again with
Lemma~2.16~b) applied to $\chi_k$, that $\vert h_{k} (z) \vert \leq
4 (m+1) \vert z \vert$ on $B\!\left( 0, \frac{r_k}{16} \right)$.
Again applying the maximum principle to $\frac{h_k}{z^{m+1}}$ on $B
\left( 0, \frac{r_k}{16} \right)$ we get Claim 1.

\medskip\noindent
We introduce several auxiliary functions in (i), (ii) and (iii)
below.

\medskip\noindent
(i) For $k \geq 0$ we define $A_k \in C^0 (T'_{k+1})$ by $A_k (z) :=
\sum\limits_{\gamma \leq R} a_{\gamma} P_{\gamma} (\log (\tau_{k}
(z))) (\tau_{k} (z))^{\gamma}$.\\
Given $k \geq 1$ we see by Lemma 2.15 that $u_k :=
\overline{\sum\limits_{\ell = 1} c_{k, \ell} A_{k}^{\ell}} \in
{\mathbb C} [[z^*]]^{\fin}_{\log}$. Note that each element of $\supp
(u_k)$ is a linear combination of the elements of $\supp (g)$ with
positive integers as coefficients. We write $u_{k} =
\sum\limits_{\gamma > 0} b_{k, \gamma} P_{k, \gamma} (\log z)
z^{\gamma} \in {\mathbb C} [[z^*]]_{\log}^{\fin}$. We set
$$
v_{k} :=  \sum\limits_{\gamma \leq R} b_{k, \gamma} P_{k, \gamma}
(\log z) z^{\gamma}, \;  w_{k} := \sum\limits_{\gamma > R} b_{k,
\gamma} P_{k, \gamma} (\log z) z^{\gamma}.$$

\noindent Using estimate (3) for the coefficients $c_{k, \ell}$ and
Lemma 2.15 for the logarithmic terms we find some $\widehat{L} > 1$
independent from $k$ such that after shrinking $S = S_R
> R$ if necessary
$$\vert w_k (z) \vert \leq \widehat{L}^k \; \vert z \vert^S \quad \mbox{for}
\quad z \in B_{{\mathbf L}} (s_k) \quad \mbox{and all} \quad k \geq
0. \leqno(5)$$

\noindent (ii) For $k \geq 0$ we define $Q_{k} (x, y) \in {\mathbb
C} [x, y]$ by
$$
Q_{k} (x, y)  = \sum\limits_{\ell = 1}^m c_{k, \ell} \left(
\sum\limits_{j = 1}^{\ell} \left( {\ell \atop j} \right) x^j
y^{\ell-j} \right) = \sum\limits_{\ell = 1}^m c_{k, \ell} ((x +
y)^{\ell} - x^{\ell}).$$

\noindent (iii) For $k \geq 0$ we define functions $\Omega_k$ on
$T_k \cap B_{{\mathbf L}} (s_k)$ as follows:\\
$k = 0$: We set $\Omega_0 \equiv 0$.\\
$k \to k+1$: We define
$$\Omega_{k+1} \colon T_{k+1} \cap B_{{\mathbf L}} (s_{k+1}), \; z \longmapsto
\left\{ \begin{array}{lll} \Omega_k (z) & & z \in T_k, \\
                          & \mbox{if} & \\
     Q_k (\varepsilon_k (\tau_k (z)), A_k (z))&  & z \in T'_{k+1}. \\
 \end{array} \right.$$

\bigskip
\noindent {\bf Claim 2:} For every $k \geq 0$
$$\varepsilon_k (z) = o (\vert z \vert^R) \quad \mbox{as} \quad z
\mapsto 0 \quad \mbox{on} \quad T_k \cap B_{{\mathbf L}} (s_k)
\leqno(\ast_k)$$

\noindent and
$$\vert \varepsilon_k (z) \vert \leq \vert \Omega_k
(z) \vert + 4 (m+1) \left( \frac{16}{r_k} \right)^m \vert z
\vert^{\frac{\alpha}{2} (m+1)} + \widehat{L}^k \vert z \vert^S \quad
\mbox{for} \quad z \in T_k \cap B_{{\mathbf L}} (s_k). \leqno(\ast
\ast_k)$$

\noindent We prove Claim 2 by induction on $k$.

\noindent $k = 0$: The base case is obvious by (2) and the
definition of $\Omega_0$.

\noindent $k \to k + 1$: Since $\Phi_{k+1}$ extends $\Phi_k$ and
since $\Omega_{k+1}$ extends $\Omega_k$ the claim holds for $z \in
T_k \cap B_{{\mathbf L}} (s_{k+1})$. Let $z \in T'_{k+1} \cap
B_{{\mathbf L}} (s_{k+1})$. Then

$$\begin{array}{l}
\Phi_{k+1} (z) = \overline{\chi_{k} (\Phi_{k} (\tau_k (z)))} \\
              \\
= \overline{\sum\limits_{\ell = 1}^m c_{k, \ell} \left(A_k (z) +
\varepsilon_{k} (\tau_k (z)) \right)^{\ell} + h_{k} (\Phi_{k} (\tau_k (z)))} \\
                 \\
= \overline{\sum\limits_{\ell = 1}^m c_{k, \ell} \left( A_k (z)
\right)^{\ell}} + \overline{Q_{k} \left( \varepsilon_{k} (\tau_k
(z)), \; A_k
(z) \right)} + \overline{ h_{k} (\Phi_{k} (\tau_k (z)))} \\
                      \\
\hspace{-2.7cm}(6) \hspace{2.1cm}= v_{k} (z) + w_{k} (z) +
\overline{\Omega_{k+1} (z)} + \overline{h_{k} (\Phi_{k} (\tau_k (z)))}.\\
\end{array}$$

\noindent To show $(\ast_{k+1})$ and $(\ast\ast_{k+1})$ we prove the
following

\bigskip
\noindent {\bf Claim 3:} Assuming that $(\ast)_k$ holds we get that
$\Phi_{k+1} (z) - v_k (z) = o ( \vert z \vert^R)$ as $z \mapsto 0$
on $T'_{k+1} \cap B_{{\mathbf L}} (s_{k+1})$.

\bigskip\noindent
Proof of Claim 3:\\
Since $S > R$ we get that $w_k (z) = o ( \vert z \vert^R)$ by (5).
Since $m \frac{\alpha}{2} > R$ we get that $h_k (\Phi_k (\tau_k
(z))) = o ( \vert z \vert^R)$ by (1) and (4). Applying $(\ast_k)$ to
$\varepsilon_k$ we deduce by the definition of $Q_k$ and $A_k$ (see
also the remarks following Fact 2.4) that $\Omega_{k+1} (z) = o (
\vert z \vert^R)$. This proves Claim 3 by (6).

\medskip\noindent We continue with the induction step of Claim 2.
For $z \in T_k \cap T'_{k+1}$ we have $z = \tau_k (z)$ and
$\Phi_{k+1} (z) = \Phi_k (z)$. Therefore we conclude by applying
Claim 3, using the inductive hypothesis $(\ast_k)$, and by the
uniqueness of the asymptotic expansion that $v_k =
\sum\limits_{\gamma \leq R} a_{\gamma} P_{\gamma} (\log z)
z^{\gamma}$. Hence $\varepsilon_{k+1} (z) = \Phi_{k+1} (z) - v_k (z)
= w_k (z) + \overline{\Omega_{k+1} (z)} + \overline{h_k (\Phi_k
(\tau_k (z)))}$ for $z \in T'_{k+1} \cap B_{{\mathbf L}} (s_{k+1})$.
The second equality holds by (6). The first equality and Claim 3
give $(\ast_{k+1})$. We obtain $(\ast\ast_{k+1})$ by the second
equality and by (1), (4) and (5). So Claim 2 is proven.

\medskip\noindent
By the definition of $r_k$ and by $(\ast\ast_k)$ we find some $L >
1$ independent from $k$ such that
$$
\vert \varepsilon_k (z) \vert  \leq \big\vert \Omega_k (z) \big\vert
+ L^k \vert z \vert^T \quad \mbox{for} \quad z \in T_k \cap
B_{{\mathbf L}} (s_k) \leqno(7)$$

\noindent where $T : = \min \left\{ \frac{\alpha}{2} (m+1), S
\right\}$. Note that $T > R$. Using estimate (3) for the
coefficients $c_{k, \ell}$ and Lemma~2.15 for the logarithmic terms
we can enlarge $L$ (independently from $k$) such that for all $k
\geq 1$ and $z \in T'_k \cap B_{{\mathbf L}} (s_k)$
$$\begin{array}{l}
\big\vert \Omega_k (z) \big\vert = \big\vert Q_{k-1} \left(
\varepsilon_{k-1} (\tau_{k-1} (z)), \sum\limits_{\gamma \leq R}
a_{\gamma} P_{\gamma} (\log (\tau_{k-1} (z)))
(\tau_{k-1} (z))^{\gamma} \right) \big\vert \\
\leq L^k (\vert \varepsilon_{k-1} (\tau_{k-1} (z)) \vert + \vert
\varepsilon_{k-1} (\tau_{k-1} (z)) \vert^m).
\end{array} \leqno(8)$$

\noindent We define $D_0 : = L$ and recursively $D_{k+1} : = 3 L^k
D_k$. By induction on $k \geq 0$ we show that $\vert \varepsilon_{k}
(z) \vert \leq D_k \vert z \vert^T$ on $T_k \cap B_{{\mathbf L}}
(p_k)$ where $p_0 := s_0$ and $p_k : = \min \{ s_k, D_{k-1}^{-T} \}$
for $k \geq 1$.

\medskip\noindent
$k = 0$: The base case is a consequence of (2).

\noindent $k \to k+1$: By the inductive hypothesis we have $\vert
\varepsilon_k (\tau_k (z)) \vert \leq D_k \vert z \vert^T$ on
$T'_{k+1} \cap B_{{\mathbf L}} (p_k)$. By the definition of $p_k$
and $D_k$ we see that $\vert \varepsilon_k (\tau_k (z)) \vert \leq
1$ on $T'_{k+1} \cap B_{{\mathbf L}} (p_k)$. Using this and (8) we
obtain that $\big\vert \Omega_{k+1} (z) \big\vert \leq 2 L^k D_k
\vert z \vert^T$. With (7) we obtain $\vert \varepsilon_{k+1} (z)
\vert \leq 3 L^k D_k \vert z \vert^T$ and are done.

\medskip\noindent
By the definition of $D_k$ and $p_k$ we find some $M > 1$ such that
$\vert \varepsilon_{k} (z) \vert \leq M^{k^2} \vert z \vert^T$ on
$T_k \cap B_{{\mathbf L}} (q_k)$ where $q_k : = M^{-k^2} \leq p_k$,
$k \geq 0$. We choose $R < \overline{T} < T$ with $T - \overline{T}
\leq 1$. We set $\overline{q}_k : = M^{-
\frac{k^2}{T-\overline{T}}}$ and obtain on $T_k \cap B_{{\mathbf L}}
(\overline{q}_k)$
$$
\big\vert \Phi_{k} (z) - \sum\limits_{\gamma \leq R} a_{\gamma}
P_{\gamma} (\log z) z^{\gamma} \big\vert = \vert \varepsilon_{R, k}
(z) \vert  \leq M^{k^2} \vert z \vert^{T- \overline{T}} \vert z
\vert^{\overline{T}} \leq \vert z \vert^{\overline{T}}.$$

\noindent Using a similar argument as at the end of the proof of
Theorem~2.17 we find some $\overline{M} > 1$ such that
$$\big\vert \Phi (z) - \sum\limits_{\gamma \leq R} a_{\gamma} P_{\gamma}
(\log z) z^{\gamma} \big\vert \leq \vert z \vert^{\overline{T}}$$

\noindent on the set $\{ (r, \varphi) \in {\mathbf L} \; \vert \;
\varphi
> 0$ and $r < \overline{M}^{- (\log_+ \varphi)^2} \}$ where $\log_+
\varphi : = \max \{ 1, \log \varphi \}$. Repeating the reflection
process in the negative direction we see that
$$\Phi (z) - \sum\limits_{\gamma \leq R} a_{\gamma} P_{\gamma}
(\log z) z^{\gamma} = o (\vert z \vert^R) \quad \mbox{as} \quad z
\to 0$$

\noindent on some admissible domain $U_R \subset U$ since
$\overline{T} > R$ and $(\log \varphi)^2 \leq \sqrt{\varphi}$
eventually.

\medskip\noindent
So we have proven that $\Phi \in Q^{\log}$. If $\sphericalangle
D/\pi \in {\mathbb R} \setminus {\mathbb Q}$ then the asymptotic
expansion $g$ of $\Phi$ is an element of ${\mathbb C} [[z^*]]$ by
Fact 2.4. The proof shows that $T \Phi = g$ and therefore $\Phi \in
{\mathcal Q}$. \hfill $\square$

\bigskip\noindent
{\bf Remark 2.19.}

\noindent The proof of the theorem of Lehman (see Fact 2.4 and the
subsequent remark) and the proofs of the preceeding Theorems 2.17
and 2.18 show actually the following: let $D$ be a domain with an
analytic corner, let $r
> 0$ and let $\Phi \colon {\mathbb H} \cap B (0, r) \to D$ be a
holomorphic map with the following
properties:\\
a) $\Phi$ has a continuous extension to $\overline{{\mathbb H}} \cap
B (0,
r)$ with $\Phi (0) = 0$.\\
b) $\Phi ([0, r [) \subset \Gamma_1$ and $\Phi (]-r, 0]) \subset
\Gamma_2$.

\noindent Then there is a quadratic domain $U \subset {\mathbf L}$
such that $\Phi \in {\mathcal Q}^{\log} (U)$. If $\sphericalangle
D/\pi \in {\mathbb R} \setminus {\mathbb Q}$ then $\Phi \in
{\mathcal Q} (U)$.

\bigskip\noindent
{\bf Proof of Theorem 2.12:}

\noindent The Riemann map $\Phi \colon {\mathbb H} \to \Omega$ maps
0 to a boundary point $x \in \partial \Omega$ with attached angle
$\sphericalangle$ greater than 0. After some translation of the
domain $\Omega$ we can assume that $x = 0$. Let $C$ be a simply
connected and semianalytic domain such that $C$ is representative of
the germ of $\Omega$ at 0 to which the germ of ${\mathbb H}$ at 0 is
mapped by $\Phi$. We choose $C$ such that $\partial C \cap
\partial \Omega$ consists of two semianalytic branches $\Gamma_1$ and
$\Gamma_2$ with $\Gamma_1 \cap \Gamma_2 = \{ 0 \}$. Note that
$\sphericalangle_0  C = \sphericalangle$. Moreover, we choose $C$ in
such a way that after some rotation $\rho$ there is a convergent
Puiseux series $\psi_1 \colon [0, \delta [ \to {\mathbb R}$ with
$\Gamma_1^* = \{ (t, \psi_1 (t)) \; \vert \; 0 \leq t < \delta \}$
where $\Gamma_1^* := \rho (\Gamma_1)$. This can be done by analytic
cell decomposition (see [10, pp.508-509]) and the fact that
subanalytic functions in one variable are given by such series (see
for example [7, p.192]). There is some $d \in {\mathbb N}$ and some
convergent real power series $\chi_1 \colon ] - \delta^d, \delta^d [
\to {\mathbb R}$ such that $\psi_1 (t) = \chi_1 (t^{\frac{1}{d}})$,
$t \geq 0$. Hence $\Gamma_1^* = \{ (t^d, \chi_1 (t)) \; \vert \; 0
\leq t < \delta^d \}$. We consider $\varphi^*_1 \colon B (0,
\delta^d) \rightarrow {\mathbb C}$, $z \mapsto z^d + i \chi_1 (z)$.
Then $\varphi_1^* \in {\mathcal O} (B(0, \delta^d))$ and $\Gamma_1^*
= \varphi_1^* ([ 0, \delta^d [)$. By the same argument applied to
$\Gamma_2$ we find (after back-rotation and some dilatation)
holomorphic functions $\varphi_1, \varphi_2 \in {\mathcal O} (B(0,
1))$ with $\varphi_1 (0) = \varphi_2 (0) = 0$ such that $\Gamma_1 =
\varphi_1 ([0, 1[)$ and $\Gamma_2 = \varphi_2 ([0, 1 [)$. Let $m_1,
m_2 \in {\mathbb N}$ be such that $\varphi_1 (z) = z^{m_1}
\widehat{\varphi}_1 (z)$, $\varphi_2 (z) = z^{m_2}
\widehat{\varphi}_2 (z)$ where $\widehat{\varphi}_1,
\widehat{\varphi}_2 \in {\mathcal O} (B(0, 1))$ and
$\widehat{\varphi}_1 (0) \cdot \widehat{\varphi}_2 (0) \ne 0$. We
can replace $\varphi_2 (z)$ by $\varphi_2 (z^{m_1})$ and can
therefore assume that $m_1$ divides $m_2$.\\
Let $r > 0$ such that $\Phi (\overline{{\mathbb H}} \cap B (0, r))
\subset \overline{C}$ and $\Phi ([0, r [) \subset \Gamma_1$, $\Phi
(]- r, 0]) \subset \Gamma_2$ (we switch 1 and 2 if necessary). We
apply finitely many elementary transformations to $\Phi$ and $C$. We
obtain functions $\Phi ^{(i)}$ and domains $C^{(i)}$, $i = 1, 2, 3$,
such that $C^{(3)}$ is a domain with an analytic corner. The domains
$C^{(i)}$ allow a similar description as $C$. We denote the
corresponding data describing $C^{(i)}$ by $\Gamma_{1,2}^{(i)}$,
$\varphi_{1, 2}^{(i)}$ and $m_{1,2}^{(i)}$.

\begin{itemize}
\item[1)] We consider $\Phi^{(1)} \colon {\mathbb H} \cap B(0, r) \to
C^{(1)}$, $z \mapsto \Phi (z)^{\frac{1}{m_1}}$, where $C^{(1)} : =
(C)^{\frac{1}{m_1}}$ (we take an appropriate $m_1$-th root of $\Phi$
and on $C$; note that $C$ is simply connected). Since $m_1$ divides
$m_2$ we see that $C^{(1)}$ has a similar description as $C$ but
additionally $m_1^{(1)} = 1$. Moreover, $\sphericalangle_0 C^{(1)} =
\sphericalangle/m_1$.
\item[2)] We consider (after shrinking $r > 0$ if necessary)
$\Phi^{(2)} \colon {\mathbb H} \cap B (0, r) \to C^{(2)}$, $z
\mapsto - (\varphi_1^{(1)})^{-1} \Phi^{(1)} (z)$, where $C^{(2)} : =
- (\varphi_1^{(1)})^{-1} (C^{(1)} )$. Note that $\varphi_1^{(1)}$ is
invertible at 0 since $m_1^{(1)} = 1$. We can choose $C$ a priori
such that $(\varphi_1^{(1)})^{-1} \in {\mathcal O} (B(0, s))$ and
$C^{(1)} \subset B(0, s)$ for some $s > 0$. We have that
$\Gamma_1^{(2)} \subset {\mathbb R}_{\leq 0}$. Moreover,
$\sphericalangle_0 C^{(2)} = \sphericalangle_0 C^{(1)}$.
\item[3)] We consider $\Phi^{(3)} \colon {\mathbb H} \cap B (0, r) \to
C^{(3)}$, $z \mapsto \rho (\Phi^{(2)} (z))^{\frac{1}{m_2^{(2)}}}$,
where $C^{(3)} : = \rho (C^{(2)})^{\frac{1}{m_2^{(2)}}}$ and $\rho
\in {\mathbb C}$ with $\vert \rho \vert = 1$ such that
$\Gamma_1^{(3)} \subset {\mathbb R}_{\leq 0}$. Then $m_2^{(3)} = 1$.
Moreover, $\sphericalangle_0  C^{(3)}  = \sphericalangle_0 C^{(2)}
/m_2^{(2)}$.
\end{itemize}

\noindent By construction $C^{(3)}$ has an analytic corner. By
Theorem 2.18 and the subsequent Remark~2.19 we get that $\Phi^{(3)}
\in {\mathcal Q}^{\log}$ and $\Phi^{(3)} \in {\mathcal Q}$ if
$\sphericalangle \; C^{(3)} \in \pi ({\mathbb R} \setminus {\mathbb
Q})$. We have $\Phi = (\varphi_1^{(1)} (-(\rho^{-1}
\Phi^{(3)})^{m_2^{(2)}}))^{m_1}$. Ge\-ne\-ra\-li\-zing [21,
Proposition 3.9] to ${\mathcal Q}^{\log}$ and using the fact that
${\mathcal Q}^{\log}$ (resp. ${\mathcal Q}$) is a ${\mathbb
C}$-algebra (see Remark~2.10) we get that $\Phi \in {\mathcal
Q}^{\log}$ and that $\Phi \in {\mathcal Q}$ if $\Phi^{(3)} \in
{\mathcal Q}$. We have $\Phi^{(3)} \in {\mathcal Q}$ if
$\sphericalangle_0 C^{(3)}/\pi \in {\mathbb R} \setminus {\mathbb
Q}$ and the latter is the case iff $\sphericalangle_0 C / \pi \in
{\mathbb R} \setminus {\mathbb Q}$. \hfill $\square$

\vspace{1.5cm}\noindent {\large \bf 3. Riemann maps at semianalytic
domains and o-minimality}

\bigskip\noindent
Now we are able to prove Theorem A. As mentioned in the introduction
the singular boundary points are the difficult part. We use
Theorem~2.12. In [21] it was shown that the functions of ${\mathcal
Q}$ restricted to the real line generate an o-minimal structure,
denoted by ${\mathbb R}_{{\mathcal Q}}$. We show that the Riemann
map is definable in the o-minimal structure ${\mathbb R}_{{\mathcal
Q}}$ as a two variable function. We use polar coordinates.

\bigskip\noindent
{\bf Definition 3.1} (compare with [21, Definition 3.4 \& Definition
4.3]){\bf .}

\noindent Let $\lambda \in \overline{{\mathbb H}} \setminus \{ 0
\}$. We have $B( \vert \lambda \vert, \vert \lambda \vert)) \subset
{\mathbf L}$ via the identification of ${\mathbb C} \setminus
{\mathbb R}_{\leq 0}$ with ${\mathbb R}_{> 0} \times ] - \pi, \pi [
\subset {\mathbf L}$. Let $\lambda = \vert \lambda \vert \, e^{ia}$
with $0 \leq a \leq \pi$. We identify $B(\lambda, \vert \lambda
\vert)$ with $\{ (r, \varphi) \in {\mathbf L} \; \vert \; (r,
\varphi - a) \in B (\vert \lambda \vert, \vert \lambda \vert) \}$.
We set ${\mathbf t}_{\lambda} \colon B(0, \vert \lambda \vert) \to B
(\lambda, \vert \lambda \vert)$, $z \mapsto \lambda + z$, and for
$\rho > 0$ we define ${\mathbf r}^{\rho, \lambda} \colon {\mathbf L}
\times B(0, \vert \lambda \vert) \to {\mathbf L}^2$, $(z_1, z_2)
\mapsto (z_1, w_2)$, with $w_2 : = {\mathbf m} ({\mathbf p}^{\rho}
(z_1), {\mathbf t}_{\lambda} (z_2))$ (see Definition 2.5 for the
definition of ${\mathbf p}$ and ${\mathbf m}$).

\bigskip\noindent
{\bf Remark 3.2.}

\noindent Let $U \subset {\mathbf L}^2$ be a 2-quadratic domain
(compare with [21, Definition 2.4]) and let $f \in {\mathcal Q}_2^2
(U)$ (compare with [21, Definition 5.1]). Let $\lambda \in
\overline{{\mathbb H}} \setminus \{ 0 \}$. As in [21,
Proposition~4.4 \& Proposition~5.15] we find some 1-quadratic domain
$V \subset {\mathbf L} \times B_{{\mathbf L}} (\vert \lambda \vert)$
such that ${\mathbf r}^{1, \lambda} (V) \subset U$ and the function
${\mathbf r}^{1, \lambda} f : = f \circ {\mathbf r}^{1, \lambda} \in
{\mathcal Q}_1^2 (V)$.

\medskip\noindent
Theorem A gets the following precise form:

\bigskip\noindent {\bf Theorem 3.3.}

\noindent {\it Let $\Omega \subset {\mathbb C}$ be a bounded,
semianalytic and simply connected domain. Suppose that
$\sphericalangle (\Omega, x) \subset \pi ({\mathbb R} \setminus
{\mathbb Q})$ for all} $x \in \Sing (\partial \Omega)$. {\it Let $F
\colon \Omega \to B (0, 1)$ be a Riemann map. Then $F$ is definable
in the o-minimal structure ${\mathbb R}_{{\mathcal Q}}$.}

\bigskip\noindent
Proof:

\noindent We establish a local definability result: given $x \in
\overline{\Omega}$ we show that $F \big\vert_{B(x, r)}$ is definable
in ${\mathbb R}_{{\mathcal Q}}$ for some $r > 0$. Then we use
compactness of $\overline{\Omega}$ to obtain the theorem.

\medskip
\noindent Let $x \in \overline{\Omega}$.

\medskip\noindent
Case 1: $x \in \Omega$. Let $r : = \dist (x, \partial \Omega)$. Then
$F \vert_{B (x, r)}$ is real analytic, hence $F \vert_{B (x,
\frac{r}{2})}$ is definable in ${\mathbb R}_{\an}$ which is a reduct
of ${\mathbb R}_{{\mathcal Q}}$.

\medskip\noindent
Case 2: $x \in \partial \Omega \setminus \Sing (\partial \Omega)$.
Then the germ of $\Omega$ at $x$ has one or two components (compare
with Example~1.3). Let $C_*$ be a semianalytic representative of
such a component. By the Schwarz reflection principle there is some
$r > 0$ such that $F$ has a holomorphic extension to $B(x, r)$. So
$F \vert_{C_* \cap B (x, \frac{r}{2})}$ is definable in ${\mathbb
R}_{\an}$.

\medskip\noindent
Case 3: $x \in \Sing (\partial \Omega)$. Let $C$ be a component of
the germ of $\Omega$ at $x$ and let $C_*$ be a semianalytic and
simply connected domain which is a representative of $C$. If $x
\not\in \Sing (\partial C_*)$ we can argue as in Case 2. So we
assume that $x \in \Sing (\partial C_*)$. Then $\sphericalangle_{x}
C_* \in \pi ({\mathbb R} \setminus {\mathbb Q})$ by assumption. Let
$j \colon {\mathbb H} \to B(0, 1)$ be a suitable M\"obius
transformation such that $\Phi : = F^{-1} \circ j$ maps the germ of
${\mathbb H}$ at $0$ to $C$. We show that there is some $r > 0$ such
that $\Phi \vert_{{\mathbb H} \cap B (0, r)}$ is definable in
${\mathbb R}_{{\mathcal Q}}$. By Theorem~2.12 and the assumption we
have some quadratic domain $U \subset {\mathbf L}$ such that $\Phi
\in {\mathcal Q} (U)$. We define $f \colon U \times U \to {\mathbb
C}$, $(z_1, z_2) \mapsto \Phi (z_2)$. Let $a \in [0, \pi]$. We
consider $g_a : = {\mathbf r}^{1, \lambda_a} f$ with $\lambda_a : =
e^{ia}$. By Remark 3.2 we get that $g_a \in {\mathcal Q}_1^2$. We
set $G_a : = g_a (z_1, h_a (z_2))$ with $h_a (z) : = e^{i (z+a)} -
e^{ia}$. Then $G_a \in {\mathcal Q}_1^2$ by [21, Proposition~5.10].
Hence there is some $r_a
> 0$ and some quadratic domain $U_a$ such that $G_a \in {\mathcal
Q}_1^2 (U_a \times B (0, r_a))$. We can assume that $U_a = \left\{
(r, \varphi) \in {\mathbf L} \; \vert \; 0 < r < c_a \exp (- C_a
\sqrt{\vert \varphi \vert}) \right\}$ with some positive constants
$c_a, C_a$. We define $\overline{G_a} \colon U_a \times B (0, r_a)
\to {\mathbb C}$, $(z_1, z_2) \mapsto \overline{G_a (\overline{z_1},
\overline{z_2})}$ where we set $\overline{z} : = (r, - \varphi)$ for
$z = (r, \varphi) \in L$. We also denote with $\overline{z}$ the
complex conjugate of a complex number $z$. Note that $\overline{G_a}
\in {\mathcal Q}_1^2$ (compare with [21, Proposition 7.3]). We set
${\mathcal R} G_a : = \frac{1}{2} (G_a + \overline{G_a})$ and
${\mathcal J} G_a : = \frac{1}{2i} (G_a - \overline{G_a})$. Then
${\mathcal R} G_a$, ${\mathcal J} G_a \in {\mathcal Q}_{1, 1;
\varepsilon_a}$ for some $\varepsilon_a > 0$ (compare with [21,
Section 7]). Hence ${\mathcal R} G_a$ and ${\mathcal J} G_a$ are
defined on $I_a : = [ 0, \varepsilon_a] \times [- \varepsilon_a,
\varepsilon_a]$, and ${\mathcal R} G_a \vert_{I_a}$ and ${\mathcal
J} G_a \vert_{I_a}$ are definable in ${\mathbb
R}_{{\mathcal Q}}$.\\
For $(r, \varphi) \in I_a$ we get ${\mathcal R} G_a (r, \varphi) =
\mbox{Re} \, \Phi ( r e^{i (\varphi + a)}) = \mbox{Re} \, \Phi (r
\cos (\varphi + a), r \sin (\varphi + a))$ and ${\mathcal J} G_a (r,
\varphi) = \mbox{Im} \Phi (r e^{i (\varphi + a)}) = \mbox{Im} \Phi
(r \cos (\varphi + a), r \sin (\varphi + a))$. Since the polar
coordinates are definable in ${\mathbb R}_{{\mathcal Q}}$ we find by
a compactness argument (note that $a \in [0, \pi ]$) some $r > 0$,
such that $\Phi \vert_{{\mathbb H} \cap B (0, r)}$ is definable in
${\mathbb R}_{{\mathcal Q}}$. Since the M\"obius transformation is
semialgebraic we find some $s > 0$ such that $F \vert_{C_* \cap B
(x, r)}$ is definable in ${\mathbb R}_{{\mathcal Q}}$. Doing this
argument for the finitely many connected components of the germ of
$\Omega$ at $x$ we obtain the claim. \hfill $\square$

\bigskip\noindent {\bf Corollary 3.4.}

\noindent {\it Let $\Omega \subsetneqq {\mathbb C}$ be a globally
semianalytic and simply connected domain. Suppose that
$\sphericalangle (\Omega, x) \subset \pi ({\mathbb R} \setminus
{\mathbb Q})$ for all} $x \in \Sing (\partial^{\infty} \Omega)$.
{\it Let $F \colon \Omega \to B (0, 1)$ be a Riemann map. Then $F$
is definable in the o-minimal structure ${\mathbb R}_{{\mathcal
Q}}$.}

\bigskip\noindent Proof:

\noindent We can copy the proof of Theorem 3.3. For $x = \infty$ we
work with $\Omega'$ and 0, see Definition~1.4.

\noindent \hfill $\square$

\bigskip\noindent {\bf Corollary 3.5.}

\noindent {\it Let $\Omega_1$, $\Omega_2 \subsetneqq {\mathbb C}$ be
globally semianalytic and simply connected domains with
$\sphericalangle (\Omega_1, x) \subset \pi ({\mathbb R} \setminus
{\mathbb Q})$ and $\sphericalangle (\Omega_2, x) \subset \pi
({\mathbb R} \setminus {\mathbb Q})$ for all} $x \in \Sing
(\partial^{\infty} \Omega_1)$ {\it and all} $x \in \Sing
(\partial^{\infty} \Omega_2)$. {\it Then each biholomorphic map
$\Omega_1 \to \Omega_2$ is definable in ${\mathbb R}_{{\mathcal
Q}}$.}

\bigskip\noindent
{\bf Remark 3.6.}

\noindent We consider special cases of domains.

\begin{itemize}
\item[a)] If the domain in question is a polygon $P$ a Riemann map
${\mathbb H} \to P$ is a so-called Schwarz-Christoffel map (see for
example Fischer-Lieb [16, p.208]). In the case of an rectangle the
Schwarz-Christoffel map is
given by an elliptic integral.\\
For polygons we can overcome the restriction on the angles at
singular boundary points by [20, Proposition 2]: a Riemann map
${\mathbb H} \to P$ is definable in the o-minimal structure
${\mathbb R}_{\an}^{{\mathbb R}}$ which is a reduct of ${\mathbb
R}_{{\mathcal Q}}$.
\item[b)] If the domain in question is a circular polygon $C$, then a
Riemann map ${\mathbb H} \to C$ fulfills the Schwarz differential
equation (compare with [16, Theorem VI.4.4]). The solutions to the
Schwarz differential equation are exactly quotients of independent
solutions to the hypergeometric differential equation (compare with
[16, Satz VI.5.3]). Applying Theorem 3.3 we obtain that these
functions are definable in ${\mathbb R}_{\mathcal Q}$ if
$\sphericalangle_x C \in \pi ({\mathbb R} \setminus {\mathbb Q})$
for all $x \in \Sing (\partial C)$.
\end{itemize}

\vspace{0.5cm}\noindent Tobias Kaiser\\
University of Regensburg\\
Department of Mathematics\\
Universit\"atsstr. 31\\
D-93040 Regensburg\\
Germany\\
{\footnotesize Tobias.Kaiser@mathematik.uni-regensburg.de}

\end{document}